\documentclass{amsart}

\usepackage{latexsym}
\usepackage{amsmath}
\usepackage{amsfonts}
\usepackage{amssymb}
\usepackage{theorem}
\usepackage[T1]{fontenc}
\usepackage{color}
\usepackage[numbers,sort&compress]{natbib}
\usepackage{multirow}
\usepackage{graphics}

\usepackage[%
  breaklinks=true,%              allow links to break over lines
  colorlinks=true,%              Colors the text of links and anchors.
  linkcolor=red,%                Color for normal internal links.
  anchorcolor=blue,%              Color for anchor text.
  citecolor=green,%                Color for bibliographical citations in text.
  menucolor=blue,%                Color for Acrobat menu items.
  urlcolor=blue,%                 Color for linked URLs. 
  bookmarks=true,%               Whether or not to have bookmarks.
  bookmarksopen=true,%           Show bookmarks with all the subtrees expanded.
  bookmarksopenlevel=1,%         level (\maxdimen) to which bookmarks are open.
  bookmarksnumbered=true,%       Include section numbers in bookmarks
  pdfview=FitH,%                 Sets the default PDF view for each link
  pdfstartview=FitH,%             Set the startup page view 
  pdfauthor={Joshua L. Willis},%
  pdftitle={Acceleration of generalized hypergeometric functions},% 
  pdfproducer={TeXShop},%
]{hyperref}

\newcommand{\fpq}{\ensuremath{{}_pF_q\!\left(\left.\begin{smallmatrix}\alpha_1,\ldots,\alpha_p\\ 
                  \beta_1,\ldots,\beta_q\end{smallmatrix}\right|z\right)}}
\newcommand{\Fpq}{\ensuremath{{}_pF_q\!\left(\left.\begin{matrix}\alpha_1,\ldots,\alpha_p\\ 
                  \beta_1,\ldots,\beta_q\end{matrix}\right|z\right)}}
\newcommand{\fqq}{\ensuremath{{}_{q+1}F_q\!\left(\left.\begin{smallmatrix}\alpha_1,\ldots,\alpha_{q+1}\\ 
                  \beta_1,\ldots,\beta_q\end{smallmatrix}\right|z\right)}}

\newcommand{\ftwo}{\ensuremath{{}_{2}F_1\!\left(\left.\begin{smallmatrix}\alpha_1,\alpha_{2}\\ 
                  \beta_1\end{smallmatrix}\right|z\right)}}

\numberwithin{equation}{section}

\begin{document}

\title[Acceleration of generalized hypergeometric functions]{Acceleration of generalized hypergeometric functions through
  precise remainder asymptotics}  
\author{Joshua L. Willis}
\address{Department of Physics\\
  Abilene Christian University\\
  Box 27963\\
  Abilene, Texas 79699, USA}
\address{Max-Planck-Institut f\"{u}r Gravitationsphysik\\
  Callinstr. 36\\
  D-30167 Hannover, Germany}
\email{josh.willis@acu.edu}
\keywords{generalized hypergeometric functions, series
  acceleration, recurrence asymptotics}
\subjclass[2010]{33C20, 65B10 (Primary) 33F05, 65D20 (Secondary)}

\begin{abstract}
We express the asymptotics of the remainders of the partial sums
$\{s_n\}$ of the generalized hypergeometric function $\fqq$ through an
inverse power  series $z^n n^{\lambda} \sum \frac{c_k}{n^k}$, where
the exponent  $\lambda$ and the asymptotic coefficients $\{c_k\}$ may
be recursively  computed to any desired order from the hypergeometric
parameters and  argument.  From this we derive a new series
acceleration technique that can  be applied to any such function, even
with complex parameters and at the branch point $z=1$.  For moderate
parameters (up to approximately ten) a C implementation at fixed
precision is very effective at computing these functions; for larger
parameters an implementation in higher than machine precision would be
needed.  Even for larger parameters, however, our C implementation is
able to correctly determine whether or not it has converged; and when
it converges, its estimate of its error is accurate. 
\end{abstract}

\maketitle

% MSC2010: 33C20 (Generalized hypergeometric series, ${}_pF_q$)
%          33F05 (Numerical approximation and evaluation)
%          65B10 (Summation of series)
%          65D20 (Computation of special functions, construction of tables)
%
% PACS2010: 02.30.Gp (Special functions)
%           02.30.Lt (Sequences, series, and summability)
%           02.30.Mv (Approximations and expansions)
%           02.60.Gf (Algorithms for functional approximation)

\section{Introduction}
\label{sec:introduction}

The generalized hypergeometric function $\fpq$ is ubiquitous in
applied mathematics; a wide array of special functions are
particular cases of this function.  Hence the numerical evaluation of
this function is an important problem. In many instances, specialized
methods for a particular special function are the most computationally
efficient, but there are still situations where the generalized
hypergeometric function must be evaluated for computationally
challenging choices of parameters and argument.  In this paper we
present a new algorithm that is able to evaluate many of those
challenging cases, and we describe its numerical implementation and
testing. 

\subsection{Analytic properties of the generalized hypergeometric function}
\label{sec:analyt-prop-gener}

To understand the difficulties, we briefly recall the definition and
analytic properties of the generalized hypergeometric function; for
more details see the reference site the Digital Library of
Mathematical Functions \cite{dlmf} or its print version \cite{NIST},
the monograph of Slater \cite{slater66}, or the text of Graham
\textit{et al.} \cite{GKP94}.  Such
functions are characterized by a set of upper and lower parameters
and a single argument; in the most general case the argument and
any of the parameters may be complex.  The function can be defined
through a Taylor series about the origin, where that converges:
\begin{equation}
  \label{eq:1}
  \Fpq = \sum_{k=0}^{\infty} \frac{(\alpha_1)_k \cdots
    (\alpha_p)_k}{(\beta_1)_k\cdots(\beta_q)_k} \frac{z^k}{k!}.
\end{equation}
Here $(a)_k = a(a+1)\cdots(a+k-1)$ is the Pochhammer symbol.  The function 
${}_2F_1$ is sometimes called simply the hypergeometric function or
Gauss' hypergeometric function; the \emph{generalized} functions are
then all of those with other numbers of parameters. In this paper we
will refer to any function given by~(\ref{eq:1}) as a generalized
hypergeometric function, or sometimes simply a hypergeometric function.

As a function of the complex argument $z$, the analytic
characteristics of ${}_pF_q$ depend on the relative number of upper and lower
parameters. When $p \leq q$, the function is entire and so the power
series converges everywhere.  When $p = q+1$---the case we will we
concerned with in this paper---the series converges for $z < 1$. On
the unit circle the behavior of the series (\ref{eq:1}) with $p=q+1$
depends crucially on the real part of the parameter $\sigma$, defined as: 
\begin{equation}
  \label{eq:3}
  \sigma = \sum_{k = 0}^{q+1} \alpha_k - \sum_{k=0}^{q} \beta_k.
\end{equation}
Then the series converges at $z = 1$ if $\Re{\sigma} < 0$,
and it converges elsewhere on the unit circle if $\Re{\sigma} <
1$. 

Outside of the unit circle, ${}_{q+1}F_q$ may be defined through analytic
continuation. This may be effected through the
defining differential equation satisfied by the generalized hypergeometric
function (see \cite{dlmf}); that equation is singular at $z = 1$ and the 
function has a branch point there (logarithmic if $\sigma \in \mathbb{Z}$;
algebraic otherwise).  The branch cut is conventionally taken 
along the positive real axis from $z = 1$ to infinity.

For the Gaussian hypergeometric function ${}_2F_1$ we can in fact
evaluate the function at the branch point in closed form, thanks to
Gauss' formula
\begin{equation}
  \label{eq:7}
 {}_{2}F_1\!\left(\left.\begin{matrix}\alpha_1,\alpha_{2}\\ 
                  \beta_1\end{matrix}\right|1\right)  =
  \frac{\Gamma(\beta_1)\Gamma(\beta_1-\alpha_1-\alpha_2)}{\Gamma(\beta_1-\alpha_1)\Gamma(\beta_1-\alpha_2)}
\end{equation}
that is valid whenever $\Re{(\beta_1-\alpha_1-\alpha_2)} > 0$. 
This formula will prove very useful later for testing our method. 

Finally, though we will not be concerned with the case $p > q+1$ in
this paper, we do note that in this case the radius of convergence of
the series is zero.

With this background, consider now the numerical computation of
${}_{q+1}F_{q}$. The advice given in the standard reference
\emph{Numerical Recipes} \cite{nrc3}
is to simply use the series (\ref{eq:1}) directly when $|z| \ll 1$,
and to integrate the defining differential equation of the generalized
hypergeometric function elsewhere, of course taking care not to cross
the branch cut, and avoiding the neighborhood of the branch point.
This leaves open, however, the question of how to evaluate the
function at or near the branch point, where the differential equation
is singular and the series slowly convergent.

\subsection{Summary of results}
\label{sec:summary-results}

We address that problem in this paper through a novel series
acceleration technique.  For any ${}_{q+1}F_q$, we will show that
the sequence $\{s_n\}$ of partial sums:
\begin{equation}
  \label{eq:4}
  s_n := \sum_{k=0}^{n-1} \frac{(\alpha_1)_k \cdots
    (\alpha_{q+1})_k}{(\beta_1)_k\cdots(\beta_q)_k} \frac{z^k}{k!}
\end{equation}
satisfies an asymptotic expansion of the form:
\begin{equation}
  \label{eq:5}
  s_n \sim s + \mu\, z^n n^{\lambda} \sum_{k=0}^{\infty}
  \frac{c_k}{n^k}, \qquad n\rightarrow\infty,
\end{equation}
for undetermined constants $s$ and $\mu$ but for asymptotic
coefficients $c_k$ that can be computed to any desired order
recursively from the hypergeometric parameters
$\alpha_1,\ldots,\alpha_{q+1}$, $\beta_1,\ldots,\beta_q$ and the complex
argument $z$.  The exponent $\lambda$ will be shown to be:
\begin{equation}
  \label{eq:6}
  \lambda =
  \begin{cases}
    \sigma & \text{when $z=1$} \\
    \sigma-1 & \text{otherwise}
  \end{cases}
\end{equation}
with $\sigma$ defined as in (\ref{eq:3}).   

Using the asymptotic expansion (\ref{eq:5}), we may use any two
successive computed partial sums $s_n$ and $s_{n+1}$ to estimate the
undetermined coefficients $s$ and $\mu$; the estimate for $s$ then
becomes the accelerated estimate for the sum and hence the function
$\fqq$ itself. We describe in some detail a numerical implementation
of this technique, which estimates both truncation and floating point
errors to determine either that the algorithm has converged to the
user specified tolerance, or that convergence is impossible at
standard machine precision.  Even in the latter case, however, the
algorithm itself is still capable of accelerating convergence; it just
must be implemented at higher than machine precision.  We present
tests of this algorithm to show that it is robustly able to either
accelerate convergence (in many cases dramatically) or correctly
conclude that a higher working precision is needed.  The algorithm is
particularly effective at the branch point $z=1$.

There are two key features of this work that deserve highlighting:
\begin{enumerate}
\item Our emphasis is on a robust algorithm that can handle complex
  parameters and argument.  We want it to succeed as often as
  possible---without user intervention to determine convergence---and to
  reliably indicate failure when it has not converged.  Considerable
  effort has therefore been spent in designing error estimation and
  stopping criteria, and the tests summarized in
  section~\ref{sec:testing-method} are designed to thoroughly probe
  how well the algorithm meets these criteria.  
\item The key novelty of the algorithm is its analytic calculation of
  the remainder asymptotics.  This is possible only because we have an
  analytic expansion of the term ratio in inverse powers.  Thus, this
  algorithm requires detailed analytic knowledge about the series
  accelerated.  This is a strength of the method in that we might hope
  (and in fact will see) that specific analytic knowledge about our
  series allows our method to succeed where other methods fail.  But
  it is also a limitation, since there is no obvious way to generalize
  the method to series that are only known numerically (certainly a
  very important class). Nevertheless, even within this limitation we
  believe that there is more to be explored, as many functions and
  series beyond the ${}_{q+1}F_q$ functions we consider here may be
  amenable to this approach.
\end{enumerate}

\subsection{Previous work}
\label{sec:previous-work}

The literature on computing the generalized hypergeometric function
depends on the restrictions placed on the number and values of the
parameters. For real parameters and argument to $\ftwo$, it is often
efficient to piece together different approximations based on the
values of the parameters and the argument; this is the approach taken,
for example, by the popular GNU Scientific Library \cite{GSL09}.  A
more detailed analysis of the algorithms appropriate for different
(real) parameters may be found in the work of Muller \cite{muller01}
for the particular case of the confluent hypergeometric function
${}_1F_1$.  Pearson's master's thesis \cite{pearson09} considers both
the confluent hypergeometric function and ${}_2F_1$, and moreover
considers complex parameters and argument. The paper \cite{forrey97}
of Forrey describes software that uses functional transformations and
difference equations to evaluate the Gaussian hypergeometric function
for arbitrary real argument, and Becken and Schmelcher \cite{BS00}
consider analytic continuation formulae to again customize the
computation based on the range of the argument. Chatterjee and Roy
\cite{CR08} consider a modification of standard Levin-type methods
tailored to the hypergeometric function, and Gautschi
\cite{gautschi02} considers evaluation of both Gaussian and confluent
hypergeometric functions for complex arguments, but real parameters,
using Gaussian quadrature to evaluate integral representations of the
functions. Weniger looked at using traditional series acceleration
methods but irregular input data in \cite{weniger01} and considered
divergent hypergeometric series at  $z=-1$ using a method tailored to
alternating series in \cite{weniger08}.  Finally Kalmykov
\cite{kalmykov04} and Kalmykov \textit{et al.} \cite{KWY07} consider an
expansion of the Gaussian function near integer values of its
parameters.  

For complex parameters and argument for the generalized hypergeometric
function $\fqq$, the literature is much more sparse. Skorokhodov
considers analytic continuation via symbolic manipulations for the 
generalized hypergeometric function in the neighborhood of $z=1$ in
\cite{skoro03,skoro04}. Ferreira \textit{et al.} \cite{FLS06} consider an
expansion valid for larger lower parameters of the Gaussian
function. Aside from the recommendation of \cite{nrc3} already given
above, Perger \textit{et al.} \cite{PBN93} implement the defining
series (\ref{eq:1}) directly, in higher than machine precision.  While
this can protect against some instances of floating point error, it
does not speed the convergence of the series itself.  The most
sophisticated software package to handle the generalized
hypergeometric function that the present author has encountered is
Johansson's \texttt{mpmath} Python module \cite{mpmath}.  That package
uses a mix of direct series calculation and analytic continuation, as
well as the Shanks' series acceleration method near the unit circle
but away from $z = 1$.  Near the branch point Euler-Maclaurin
summation is used.  

Series acceleration has a long history; we will give a brief review
and more pointers to the literature in 
section~\ref{sec:surv-series-accel}. Here we only mention a few recent
techniques that have been specifically applied to hypergeometric
functions. Wozny and Nowak \cite{WN09} and Wozny \cite{wozny10}
consider a new series acceleration technique that they apply
(among other examples) to some instances of the generalized
hypergeometric function.  Their approach is based on finding certain
difference operators that approximately annihilate the remainder term
of the series. Paszkowski \cite{pasz11} considers several acceleration
algorithms and how they may be modified if an asymptotic form for the
partial sums is known. While it is mentioned that the generalized
hypergeometric functions belong to this class, an explicit expression
for the asymptotic coefficients is not given, and only a few low order
expansions are considered in examples.  Likewise Lewanowicz and
Paszkowski \cite{LP95} consider an acceleration method based 
on the asymptotic expansion of the term ratio, as we will, but those
authors do not apply it to the asymptotics of the partial sum, and
they are only able to accelerate certain parameter choices for
${}_3F_2$ functions at $z=\pm 1$.  In Skorokhodov \cite{skoro05} and
Bogolubsky and Skorokhodov \cite{BS06} the authors
use an asymptotic expansion of the terms---rather
than the term ratios---to calculate an approximant to the truncation
error that they evaluate using Hurwitz zeta functions.  We shall
compare the performance of our method to theirs in section
\ref{sec:zeta-funct-accel}. 

Perhaps closest in spirit to the present work is that of Weniger
\cite{weniger07}, which is based on finding asymptotic approximations
to the remainder terms of a partial series summation through symbolic
linear algebra. That method can in some cases yield analytic
expressions, as for the Dirichlet series of the Riemann zeta function
(section 5 of \cite{weniger07}) or the divergent Euler series for the
exponential integral (as in Borghi \cite{borghi10}). However, the
method presented in \cite{weniger07} is challenging even for
${}_2F_1$, particularly as compared to the complete asymptotic series
that we will give in this paper.    

Finally, we mention that B\"uhring has considered both the behavior of
the generalized hypergeometric functions near unity \cite{buhring92},
seeking to find the analogue of Gauss' formula (\ref{eq:7}) for
higher order functions; and separately considered the asymptotic
behavior of the partial sums of generalized hypergeometric functions
at unity \cite{buhring03}.  In each case, the expressions derived are
nested infinite sums of hypergeometric functions of lower order, so
the results do not seem well adapted to numeric computation.  More
practically useful for us is his work in \cite{buhring87,BS98} and the
earlier work of N{\o}rlund \cite{norlund55}, which give expansions valid
near the branch point.  They can conceivably be leveraged to take an
efficient evaluation method at the branch point and evaluate a
hypergeometric function near the branch point using expansions valid
in a neighborhood of the branch point.

\section{Accelerating the convergence of the series}
\label{sec:rev-gener-hyperg}

The well-known Euler's method shows that series acceleration dates to
at least the eighteenth century, and both Knopp \cite{knopp64} and
Tweddle \cite{tweddle03} cite Stirling as the earliest to develop a
series acceleration method, but the last several decades have seen the
development of a variety of sophisticated techniques.  For reviews,
see the articles of Brezinski \cite{brezinski85}, Homeier
\cite{homeier00}, and Weniger \cite{weniger89}, and the
monographs of Brezinski \cite{brezinski77,brezinski78}, Brezinski and
Redivo Zaglia \cite{BR91}, Sidi \cite{sidi03}, Walz \cite{walz96} and
Wimp \cite{wimp81}.  We will give just enough background in the next
section to place our new method in context.

\subsection{Review of series acceleration methods}
\label{sec:surv-series-accel}

The basic idea of any series acceleration technique is to use the
expected form of the partial sums $\{s_n\}$ of a series---which by  
themselves may be slowly convergent or even divergent---to create
a new sequence $\{s_n'\}$ that converges to the same limit $s$ (or, in
the case of a divergent series, antilimit), but that does so more
rapidly, in the sense that:
\begin{equation}
  \label{eq:8}
  \lim_{n\rightarrow\infty} \frac{s_n'-s}{s_n-s} = 0.
\end{equation}

To motivate these transformations, we start from the explicit
expression of the sequence $\{s_n\}$ of partial sums by means of the
terms $\{t_k\}$ of the series:
\begin{equation}
  \label{eq:11}
  s_n := \sum_{k=0}^{n-1} t_k.
\end{equation}
Now partition the partial sum into
its (anti-)limit $s$ and the \emph{remainder} $\rho_n$:
\begin{equation}
  \label{eq:12}
  \rho_n := -\sum_{k=n}^{\infty} t_k
\end{equation}
so that:
\begin{equation}
  \label{eq:13}
  s_n = s + \rho_n.
\end{equation}
Various algorithms can then be devised by approximating the remainder
as $\rho_n \approx \omega_n \mu_n$, where $\omega_n$ is an explicit remainder
estimate, and $\mu_n$ is an $O(1)$ correction factor containing $m$
free parameters.  The order of the transformation is $m$.

For instance, if the terms of the series are alternating, then a
natural estimate of the remainder would be the first term not
included, $t_n$.  If we choose $\omega_n = t_n$ and choose:
\begin{equation}
  \label{eq:14}
  \mu_n = \sum_{k = 0}^{m-1} \frac{c_k}{(n+\beta)^k}
\end{equation}
for some positive $\beta$, and undetermined coefficients $c_k$, then
from any sequence of $m+1$ successive partial sums, we may determine
values for the $m$ coefficients $c_k$ and the limit $s$.  This choice
of remainder estimate $\omega_n$ and correction factor $\mu_n$ gives
rise to Smith and Ford's modification \cite{SF79} of Levin's $t$
transformation~\cite{levin73}; other choices for either $\omega_n$ or
$\mu_n$ give other sequence transformations; for details see the
review articles and monographs mentioned earlier.  The actual
estimation of $s$ from the $m+1$ partial sums can be expressed as a
ratio of determinants, but is more commonly implemented
recursively~\cite{weniger89,homeier00}. 

Among all of these methods, we will single out one known as the
$E$-method \cite{BR91,brezinski80,brezinski82,havie79}, because it
largely subsumes all of the others as special cases: starting from a
set of functions $\{g_i(n)\}_{i=1}^{m}$, one requires:
\begin{equation}
  \label{eq:2}
  s_n = s + \sum_{i=0}^{m} c_i g_i(n).
\end{equation}
From the $m+1$ sums $s_n, s_{n+1}, \ldots, s_{n+m}$ we can calculate
$s$ and the $c_i$; the value of $s$ is then the accelerated sum
determined by the $E$-method for that particular choice of functions
$g_i$ and those partial sums. We will compare the performance of this
method to that of this paper in section \ref{sec:e-method}.

The effectiveness of such algorithms depends on the nature of the
series to be summed; for instance, we mentioned above that the $t$
transformation was designed with alternating 
series in mind.  One central characteristic of convergent series that
influences the effectiveness of acceleration methods is whether that
convergence is \emph{linear} or \emph{logarithmic}.  The
former means that if the limit of the sequence of partial sums is $s$,
then
\begin{equation}
  \label{eq:15}
  \lim_{n\rightarrow\infty} \frac{s_{n+1}-s}{s_n-s} = \ell
\end{equation}
with $0 < |\ell| < 1$.  On the other hand, if $\ell = 1$, then the
convergence is logarithmic~\cite{weniger89}.

The distinction is important because while several series acceleration
methods can be shown to accelerate the convergence of \emph{any}
linearly convergent series~\cite{SF79}, it is known from the work of
Delhaye and Germain-Bonne \cite{DG82} that
no method is able to accelerate the convergence of all logarithmically
convergent series.  Yet at the branch point of the generalized
hypergeometric function (as we will show in the
next section) the convergence of the series is logarithmic.
Section 14 of \cite{weniger89} contains a summary of methods that may
be applied to logarithmically convergent series; while alternating series
are often tractable, with a generic choice of complex parameters the
terms in the generalized hypergeometric series are complex and exhibit
no particular sign pattern (see Sidi \cite{sidi06} for examples of the
effect of irregular sign patterns on series acceleration techniques).
Moreover, some series acceleration techniques are quite sensitive to
the details of the asymptotic form of the remainders \cite{weniger89};
they are able to accelerate convergence when $\rho_n \sim n^{-k}$ for
an integer $k$, but fail for non-integral but real exponents.  For
complex exponents---as we will encounter for generalized
hypergeometric sums---very little indeed seems to be known.   

Thus, while there is a wide array of sophisticated series acceleration
methods able to speed the convergence, often substantially, of
many series, and even to sum many strongly divergent series, there
does not seem to be a method that is broadly applicable to the series
expansion of the generalized hypergeometric function, particularly
at its branch point.  Even very new methods, such as those in
\cite{weniger07,WN09,wozny10}, typically consider only real
hypergeometric parameters.

Instead, we propose a method that can be applied to any generalized
hypergeometric of the form ${}_{q+1}F_q$.  Unlike conventional series
acceleration methods, instead of a simple choice for the remainder
estimate $\omega_n$ and a complicated choice for the correction
$\mu_n$ (so that several successive partial sums are needed to
calculate each $s_n'$) we will use a very precise remainder estimate
but just a constant for our correction factor. Thus, we will need
only two successive partial sums to compute each estimate of the
series limit.  We turn now to the  derivation of that method. 

\subsection{Derivation of the partial sum asymptotics}
\label{sec:deriv-partial-sum}

By comparing the general form (\ref{eq:11}) to the expression
(\ref{eq:4}) for the partial sums of the generalized 
hypergeometric and making use of the definition of the
Pochhammer symbol, we may write the ratio of two successive terms of
the generalized hypergeometric function as:
\begin{equation}
  \label{eq:16}
  \frac{t_{k+1}}{t_k} = z
  \frac{(\alpha_1+k) \cdots (\alpha_{q+1}+k)}{(\beta_1+k) \cdots
    (\beta_q+k)(1+k)} := z\,r(k);
\end{equation}
that is, the ratio of two successive terms is a \emph{rational}
function of the index $k$.  This property of generalized
hypergeometric functions is very well known; indeed, any function
whose term ratio is a rational function of the term index may be
expressed in terms of generalized hypergeometric functions
\cite{GKP94}. 

Because the ratio of terms satisfies a \emph{first}-order recurrence
relation, both the remainders $\{\rho_n\}$ and the partial sums
$\{s_n\}$ satisfy a \emph{second}-order recurrence; in fact, the two
sequences satisfy the \emph{same} second order recurrence.  In equations, we
have:
\begin{equation}
  \label{eq:9}
  \frac{t_{n+1}}{t_n} = z\,r(n) = \frac{s_{n+2}-s_{n+1}}{s_{n+1}-s_n}
  = \frac{\rho_{n+2}-\rho_{n+1}}{\rho_{n+1}-\rho_n}
\end{equation}
and therefore:
\begin{align}
  s_{n+2} -\Bigl(1+z\,r(n)\Bigr)s_{n+1} + z\,r(n) s_n &=
  0   \label{eq:10a} \\ 
  \rho_{n+2} -\Bigl(1+z\,r(n)\Bigr)\rho_{n+1} + z\,r(n) \rho_n &=
  0   \label{eq:10b} 
\end{align}

The key point is that the remainders and partial sums each satisfy a
linear, homogeneous, second-order difference equation, and the
asymptotic solutions of such equations are known very precisely.  In
principle, we could use either of equation~(\ref{eq:10a})
or~(\ref{eq:10b}) as our starting point, and typically in series
acceleration it would be more natural to focus on the behavior of the
remainders.  However, it is slightly more convenient to use
(\ref{eq:10a}) and base our results directly on the asymptotics of the
partial sums themselves.  That is because
there are two linearly independent solutions to our difference
equation, and we will find that one of those two solutions 
is always a constant, and the other approaches zero as
$n\rightarrow\infty$ when $|z| < 1$, but diverges as
$n\rightarrow\infty$ when $|z| > 1$. Were we to focus solely on the
remainders, all we could say was that when $|z| < 1$ only the
decreasing solution can be present, since we know the remainders go to
zero in that case.  By directing our attention instead to the partial
sums, we see that there is always a constant term (the value of the
function we are trying to find) as well as a remainder (the second
solution) that diverges when we are outside the radius of
convergence, and goes to zero inside the radius of convergence.  Thus,
we show directly that our asymptotic acceleration not only speeds the
convergence of the series when it does converge, but also
slows the divergence when it does not.

So consider the general problem of the asymptotic solutions (valid for
large $n$) of a linear, homogeneous, second-order difference equation
that may be written in the form:
\begin{equation}
  \label{eq:10}
  w_{n+2} + a(n) w_{n+1} + b(n) w_n = 0
\end{equation}
where the coefficient functions $a(n)$ and $b(n)$ are themselves
known, and have asymptotic expansions:
\begin{align}
  a(n) \sim \sum_{k=0}^{\infty} \frac{a_k}{n^k}, \qquad  b(n) \sim
  \sum_{k=0}^{\infty} \frac{b_k}{n^k}, \qquad n\rightarrow\infty.  \label{eq:17a} 
\end{align}
This problem has been considered by several authors, beginning with
Birkhoff \cite{birkhoff11,birkhoff30}, Adams \cite{adams28}, and
Birkhoff and Trjintzinsky \cite{BT32}.  More recently, this body of
work has been reviewed and summarized by both Wimp and Zeilberger
\cite{WZ85} and Wong and Li \cite{WL92a}, who all find Birkoff's work
notable for its complexity. 
For us, by far the most useful reference will be \cite{WL92a}, since
the authors carefully analyze the several possible cases that may
arise when finding asymptotic solutions to (\ref{eq:10}), and also
provide explicit, recursive formulas for arbitrary asymptotic
coefficients. That will be crucial.  We follow the terminology of
\cite{WL92a}, but not in general the notation.

The analysis of \cite{WL92a} only considers the case in which $b_0$
of (\ref{eq:17a}) is not zero.  As we will soon see, that will limit
the ${}_pF_q$ that our acceleration can handle to those for which $p =
q+1$.  In a later paper \cite{WL92b} the same authors consider the
more general case, and so likewise we will consider the analysis of
generalized hypergeometrics where $p \neq q+1$ in a separate work.

With that restriction on $b_0$, the authors of \cite{WL92a} show that
there are two linearly independent asymptotic solutions of
(\ref{eq:10}), which fall into one of three cases (depending on the
lowest few asymptotic coefficients of $a(n)$ and $b(n)$) as follows:
\begin{description}
\item[The normal case.]  When the two roots $\xi_1$ and $\xi_2$ to the
  \emph{characteristic equation}
  \begin{equation}
    \label{eq:17}
    \xi^2 + a_0\,\xi + b_0 = 0
  \end{equation}
  are distinct, then the two linearly independent solutions
  to~(\ref{eq:10}) are each of the form:
  \begin{equation}
    \omega_{n} \sim \xi^n n^{\lambda} \sum_{k=0}^{\infty}
    \frac{c_{k}}{n^k}, \qquad n\rightarrow\infty, \label{eq:18a}
  \end{equation}
  for the two values of $\xi$, and the exponent $\lambda$ for
  each solution depends on $\xi$ through:
  \begin{equation}
    \label{eq:19}
    \lambda = \frac{a_1\,\xi +
      b_1}{a_0\,\xi + 2b_0}. 
  \end{equation}
\item[The subnormal case.]  When the roots of the characteristic
  equation do coincide, but the double root is not the zero of the
  \emph{auxiliary equation} $a_1\,\xi + b_1 = 0$, then the two
  solutions are of the form: 
  \begin{equation}
    \label{eq:20}
    \omega_{\pm,n} \sim \xi^n e^{\pm\gamma\sqrt{n}} n^{\lambda}
    \sum_{k=0}^{\infty} \frac{c_{\pm,k}}{n^{k/2}}, \qquad
    n\rightarrow\infty,
  \end{equation}
  where $\gamma$ and $\lambda$ may be explicitly determined, but we
  will not need them.
\item[The exceptional case.] When the roots of the characteristic
  equation coincide and the double root is also the root of the
  auxiliary equation, then the two linearly independent solutions are
  again given by equation (\ref{eq:18a}), but now
  the two values of $\lambda$ are given not by (\ref{eq:19}), but
  rather as the two roots of the \emph{indicial equation}:
  \begin{equation}
    \label{eq:21}
    \lambda(\lambda-1)\xi^2 +(a_1\,\lambda+a_2)\xi+b_2 = 0.
  \end{equation}
  There are some further complications considered in \cite{WL92a} when
  the two roots of this equation either coincide or differ by a
  positive integer, but we will not need those subtleties.
\end{description}

In each of these three cases, the leading coefficient $c_0$ of the
asymptotic expansion may be taken, without loss of generality, to be
one, and the higher coefficients are then determined recursively from
the $\{a_n\}$ and $\{b_n\}$ through formulas that we will quote from
\cite{WL92a} later.  We will not discuss the derivation of these cases
and the corresponding formulas, except to say that superficially the
method is much like the series solution of differential equations: one
proposes a form of the solution and inserts it into the equation, and
this recursively determines all of the higher coefficients.  But
showing that the resulting formal solutions are indeed asymptotic is
far from trivial, and for details the reader is referred to
\cite{WL92a} and the rest of the literature cited.

We see immediately that the recursion (\ref{eq:10a}) satisfied by the
partial sums $s_n$ is of the form (\ref{eq:10}), provided we take:
\begin{equation}
  \label{eq:22}
  a(n) := -\Bigl(1+z\,r(n)\Bigr) \qquad b(n) := z\,r(n).
\end{equation}
To apply the results of \cite{WL92a}  we
need an asymptotic expansion for $r(n)$.  But for large $n$, that is
easy; because $p = q+1$ we see from (\ref{eq:16}) that $r(n)$ is a
rational function whose numerator degree equals its denominator
degree, and so we divide both numerator and denominator by $n^{q+1}$
and write:
\begin{equation}
  \label{eq:24}
r(x) :=  \frac{(1+\alpha_1x) \cdots (1+\alpha_{q+1}x)}{(1+\beta_1x) \cdots
    (1+\beta_qx)(1+x)}   
\end{equation}
where we have defined $x$ as $1/n$.  This rational function (we
deliberately use the same symbol) has a convergent Taylor series
expansion in a neighborhood of zero, and the coefficients of that
Taylor series will coincide with the asymptotic coefficients of $r(n)$
as $n\rightarrow\infty$.  We will need arbitrarily many of these
coefficients for our full acceleration method, and efficiently
calculating those is not trivial, so we defer it to the next
subsection.  However, to determine which of the three cases above we
fall under, we need only the lowest two, and elementary
calculus yields:
\begin{align}
   r_0 &= 1  \label{eq:26a} \\
   r_1 &= \sum_{k=1}^{q+1} \alpha_k - \sum_{k=0}^{q} \beta_k -1 = \sigma-1. \label{eq:26b}
\end{align}
Finally, we will also need:
\begin{align}
  a_n &=
  \begin{cases}
    -1-z& \text{if $n=0$} \\
    -z\,r_n& \text{otherwise}
  \end{cases}   \label{eq:27a} \\
  b_n &= z\,r_n  \label{eq:27b}
\end{align}
and the consequent identity $a_n = -b_n$ whenever $n \geq 1$.  Note that
$b_0 \neq 0$ precisely because we assume $p=q+1$.

Having considered the generalities, we now turn to precise formulas
for the asymptotics of the $\{s_n\}$.  As those depend on whether or
not we are at the branch point $z=1$, we take up those two cases in
turn.  Before beginning that discussion, though, we point out that
none of the results of the next two subsections apply when one of the
upper or lower parameters is a non-positive integer. 
That is because in those situations the recursion~(\ref{eq:10a}) does not
really describe the asymptotics of the partial sums as
$n\rightarrow\infty$; rather, in the first situation the hypergeometric is
a terminating polynomial, and in the second, it is undefined.  So
although our results will not apply, either case is easy
to identify and handle without series acceleration techniques.

\subsubsection{Partial sum asymptotics away from the branch point}
\label{sec:part-sum-asympt}

For our recurrence (\ref{eq:10a}), we have from
(\ref{eq:26a}--\ref{eq:27b}) that the characteristic equation
is:
\begin{equation}
  \label{eq:23}
  \xi^2 - (1+z) \xi + z = 0,
\end{equation}
and the roots of this equation are $1$ and $z$.  Thus, provided we are
not at the branch point of the hypergeometric function, the recurrence
equation is in the normal case of the three listed above.  We can also
find the corresponding exponents from (\ref{eq:19}), and we easily
determine that the exponent corresponding to $\xi = 1$ in fact
vanishes, whereas the exponent corresponding to the root $\xi=z$ is
$r_1 = \sigma-1$.

To find the asymptotic coefficients for each of these cases, we need
the recursion that those coefficients satisfy. That is given in
equation~(2.3) of \cite{WL92a}, and in our notation is:
\begin{equation}
  \label{eq:25}
  \sum_{j=0}^{k-1} \left[ \xi^2 2^{k-j} \binom{\lambda-j}{k-j} + \xi
      \sum_{i=j}^{k} \binom{\lambda-j}{i-j} a_{k-i} + b_{k-j}\right] c_j = 0.
\end{equation}
It is convenient to rewrite this equation so that it explicitly gives
$c_k$ in terms of $\{c_0,\ldots,c_{k-1}\}$.  Rearranging terms, making
use of equations (\ref{eq:27a}) and (\ref{eq:27b}), and renaming
the index of summation yields 
\begin{equation}
  \label{eq:26}
  c_k = \frac{1}{k(1-z)} \sum_{j=0}^{k-1} \left[
    2^{k+1-j}\binom{-j}{k+1-j} + \sum_{i=j}^{k+1} \binom{-j}{i-j}
    a_{k+1-i} + b_{k+1-j}\right] c_j,
\end{equation}
for the root $\xi = 1$ (and therefore $\lambda=0$), and
\begin{multline}
  c_k = -\frac{1}{k(1-z)} \sum_{j=0}^{k-1} \Biggl[
    \Bigl((2^{k+1-j}-1)z-1\Bigr)\binom{\sigma-1-j}{k+1-j}  \\ 
  \left. - z \sum_{i=j}^{k} \binom{\sigma-1-j}{i-j} r_{k+1-i} + 
    r_{k+1-j}\right] c_j
  \label{eq:27}
\end{multline}
for the root $\xi = z$ (and $\lambda = \sigma-1$). 
Each of these equations is valid for all $k \geq 1$.

In equation~(\ref{eq:27}) we have expressed the recursion in terms of
the asymptotic coefficients $r_n$ of $r(n)$, rather than the $a_n$ and
$b_n$, but we have not troubled to do that for the asymptotic
coefficients for the $\xi=1$ root.  That is because all of the
coefficients $c_k$ in equation~(\ref{eq:26}) are zero 
when $k \geq 1$.  To see this, first consider the coefficient of $c_0$
on the right hand side of~(\ref{eq:26}). The binomial coefficients
$\tbinom{0}{l}$ for any integer $l$ vanish, except for $\tbinom{0}{0}$
which equals one.  Thus, since $j < k+1$, the first term in the
square brackets vanishes because the binomial coefficient does, as do
all of the terms in the inner sum over $i$, except for the $i=j=0$
term. That term survives to give $a_{k+1}$, but that in turn cancels
the term $b_{k+1}$.  Hence, the entire coefficient of $c_0$ vanishes.

But by the recursion~(\ref{eq:26}), that means that $c_1$ vanishes,
and then by induction that all of the $c_k$ vanish when $k>0$, for the
$\xi=1$ root of the characteristic equation.  Hence, comparing
to~(\ref{eq:18a}), we see that the first asymptotic solution to the
difference equation~(\ref{eq:10a}) is simply a constant.  That is not
true for the second solution, however, and so since the partial sum
$s_n$ is a linear combination of these two asymptotic solutions, we
have shown:
\begin{equation}
  \label{eq:28}
  s_n \sim s + \mu\, z^n\,n^{\sigma -1}
  \sum_{k=0}^{\infty} \frac{c_k}{n^k}, \qquad n\rightarrow\infty,
\end{equation}
with $\sigma$ given by~(\ref{eq:3}) and the $c_k$ determined
recursively from $c_0 = 1$ by equation~(\ref{eq:27}).

This result holds so long as $z\neq 1$; in particular, it holds
outside the radius of convergence $|z|=1$.  In that case, however,
(\ref{eq:28}) shows that the remainder term diverges as
$n\rightarrow\infty$; the series is not convergent there.  It clearly
converges whenever $|z| <1$, and on the circle of convergence the
remainder is a decreasing function of $n$ provided the real part of
the exponent of $n$ is negative; that is, provided $\Re{\sigma} < 1$.
So, our asymptotic formula correctly reproduces the analytic
properties of ${}_{q+1}F_q$ whenever $z \neq 1$; we next consider $z =
1$. 

\subsubsection{Partial sum asymptotics at the branch point}
\label{sec:part-sum-asympt-1}

When $z=1$ we are no longer in the normal case of \cite{WL92a}.  To
decide between the subnormal and exceptional cases, we must examine
the auxiliary equation $a_1 \xi + b_1 = 0$.  As $a_1 = -b_1$ and $\xi
= z= 1$, our double root is in fact a root of the auxiliary equation,
so we are in the exceptional case. The indicial equation when $\xi =
1$ is: 
\begin{equation}
  \label{eq:29}
  \lambda(\lambda-1) + (a_1 \lambda +a_2) + b_2 = 0.
\end{equation}
Since $a_2+b_2 = 0$, we get:
\begin{equation}
  \label{eq:30}
  \lambda = 0 \qquad \text{or} \qquad \lambda = r_1+1 = \sigma
\end{equation}
as our two possible exponents.

Just as for the normal case, we must now examine the recursive
equations that determine the asymptotic coefficients $c_n$ in
terms of the $\{r_n\}$ for each of these two possible values of
$\lambda$. When we do so we again find that only the leading,
constant term survives for $\lambda = 0$, while for $\lambda = \sigma$
the entire series is nontrivial.  Similar considerations allow us to
avoid some special sub-cases of the exceptional case alluded to above.
Generically, the solutions to the recursion equation can have
terms logarithmic in $n$ when the difference between the two roots of
the indicial equation are an integer, but for all of the cases where
the hypergeometric series at the branch point converges (that is,
where $\Re{\sigma} < 0$) we can show through inductive
arguments  similar to that above that the logarithmic term vanishes.
Thus, the most general asymptotic solution is of the form:   
\begin{equation}
  \label{eq:31}
  s_n \sim s + \mu\,n^{\sigma} \sum_{k=0}^{\infty} \frac{c_k}{n^k},
  \qquad n\rightarrow\infty,
\end{equation}
where again $c_0 = 1$, and the higher coefficients can be shown from
equation~(7.2) of \cite{WL92a} to satisfy:
\begin{multline}
  \label{eq:32}
  c_k = \frac{1}{k(\sigma-k)} \sum_{j=0}^{k-1} \Biggl[
    \Bigl(2^{k+2-j}-2\Bigr)\binom{\sigma-j}{k+2-j} \\ 
  \left. - \sum_{i=j}^{k+1} \binom{\sigma-j}{i-j} r_{k+2-i} + 
    r_{k+2-j}\right] c_j. 
\end{multline}

As with regular points of the function, we see that the asymptotic
expansion at the branch point reproduces the correct analytic behavior
of the ${}_{q+1}F_q$ function.  Specifically, the remainder is a
decreasing function precisely when $\Re{\sigma} < 0$, the
condition we saw in section~\ref{sec:analyt-prop-gener} was needed to
ensure that the series converged at $z=1$.  We also see explicitly
from~(\ref{eq:31}) that the convergence of the series at $z=1$ is
logarithmic, with a complex critical exponent (in general).

\subsection{Recursively computing the asymptotic coefficients}
\label{sec:recurs-comp-asympt}  

The two asymptotic expansions (\ref{eq:28}) and (\ref{eq:31}),
together with the respective recursions~(\ref{eq:27})
and~(\ref{eq:32}), are a complete solution for the asymptotics of the
partial sums of the generalized hypergeometric function.  To be
numerically useful, however, we must be able to calculate arbitrary
asymptotic coefficients $r_k$ of the rational term ratio function
$r(n)$. It is here that we can make use of the explicit analytic
knowledge we have of our series; it comes to us not just as a
numerical sequence of terms. As we noted, the $r_k$ are the same as the
Taylor coefficients of the expansion around zero of the rational
function $r(x)$ of equation~(\ref{eq:24}), but that in and of itself
is not a practical solution, if we have no better means to calculate
the Taylor coefficients than by evaluating high order derivatives at
zero. 

Fortunately, there is a much more effective procedure. Define the
functions $P(x)$ and $Q(x)$ as:
\begin{equation}
  \label{eq:33}
  P(x) = (1+\alpha_1x) \cdots (1+\alpha_{q+1}x); \qquad Q(x) =
  \frac{1}{(1+\beta_1x) \cdots (1+\beta_qx)(1+x)}   
\end{equation}
so that $r(x) = P(x)Q(x)$. It is obvious that the Taylor coefficients
of $r(x)$ satisfy
\begin{equation}
  \label{eq:34}
  r_k = \sum_{j=0}^{k} P_jQ_{k-j}
\end{equation}
when $\{P_k\}$ and $\{Q_k\}$ are the Taylor coefficients of the
respective functions.

The coefficients $P_k$ and $Q_k$ can be calculated effectively because
we have already factored the functions $P(x)$ and $Q(x)$.
Specifically, rearranging some results of \cite{macd79} gives:
\begin{equation}
  \label{eq:35}
  P(x) = \sum_{k=0}^{q+1} e_k(\alpha_1,\ldots,\alpha_{q+1}) x^k
\end{equation}
and:
\begin{equation}
  \label{eq:36}
  Q(x) = \sum_{k=0}^{\infty} (-1)^k h_k(1,\beta_1,\ldots,\beta_q) x^k.
\end{equation}
In these equations, $e_k$ and $h_k$ are the \emph{elementary symmetric
  polynomials} and the \emph{complete homogeneous symmetric
  polynomials}, respectively.  These may be defined in several ways,
but the most computationally useful method defines them
recursively from the \emph{power sums} of the upper and lower
parameters (where the lower parameters must be augmented with the
implicit parameter one that is part of the definition of the
generalized hypergeometric function).  The $k$-th power sum
$p_k(x_1,\ldots,x_n)$ of any set $\{x_i\}$ of $n$ numbers is:
\begin{equation}
  \label{eq:37}
  p_k(x_1,\ldots,x_n) = \sum_{i=1}^n x_i^k.
\end{equation}
In terms of the power sums, the elementary symmetric polynomials
satisfy the recurrence:
\begin{equation}
  \label{eq:38}
  ke_k(\alpha_1,\ldots,\alpha_{q+1}) = \sum_{i=1}^{k} (-1)^{i-1}
  e_{k-i}(\alpha_1,\ldots,\alpha_{q+1}) p_i(\alpha_1,\ldots,\alpha_{q+1})
\end{equation}
and the complete homogeneous symmetric polynomials satisfy the recurrence:
\begin{equation}
  \label{eq:39}
  kh_k(1,\beta_1,\ldots,\beta_{q}) = \sum_{i=1}^{k} 
  h_{k-i}(1,\beta_1,\ldots,\beta_{q}) p_i(1,\beta_1,\ldots,\beta_{q}).  
\end{equation}
These recurrences do not determine $e_0$ or $h_0$; both are unity.

It is now straightforward to calculate the Taylor coefficients $r_k$
for any desired $k$.  We first calculate the power sums $p_i$ of the upper
and augmented lower parameters for $i$ from zero to $k$.  Then we
use the recurrence (\ref{eq:38}) to find the elementary symmetric
polynomials of the upper parameters, and the recurrence (\ref{eq:39})
to find the complete homogeneous symmetric polynomials of the
augmented lower parameters.  Finally, combining equations
(\ref{eq:34}--\ref{eq:36}) gives:
\begin{equation}
  \label{eq:40}
  r_k = \sum_{i=0}^{k} (-1)^{k-i} e_i(\alpha_1,\ldots,\alpha_{q+1})
  h_{k-i}(1,\beta_1,\ldots,\beta_q). 
\end{equation}

\subsection{Stability of the recursions}
\label{sec:stability-recursions}

The expressions for the asymptotic coefficients \eqref{eq:27}
and \eqref{eq:32}, together with the formulas of the last subsection
for calculating the Taylor coefficients, provide a determination of
the asymptotic coefficients through several nested recursions.  It is
therefore important to study the stability of these recursions.  Given
the nesting of the recursions, an analytic study is daunting, so we
investigate the stability numerically.  Even though most of our
numeric results on the performance of our algorithm are presented in
the next section, we digress briefly to study the stability of the
coefficient calculation here.  That is both because this study was
quite different from the analysis of section~\ref{sec:impl-results},
as it was conducted in much higher than machine precision, and also
because we need an assurance of the stability of the asymptotic
coefficient calculation before we can examine the overall method.

Accordingly, we used the \textbf{mpmath} \cite{mpmath} package already
mentioned in the introduction to implement the calculation of the
coefficients in fixed but arbitrary precision.  We compared the
calculation at a precision of 53 bits (standard double precision, as
we will use in the C implementation discussed in the next section) to
that calculated with twice the precision, at 106 bits.  We expect that
the accuracy of the computation will depend on the size of the
parameters, the value of $z$, and the number of coefficients
calculated, since error presumably accumulates throughout the
recursion.  To investigate these effects, we studied both $z=1$, and
$z$ chosen uniformly at random in the unit disk.  We used four upper
and three lower parameters, with real and imaginary parts each chosen
uniformly in the range $(-R,R)$, for $R \in \{1,5,10,50,100\}$.
Unlike the actual implementation of the full algorithm itself that we
study in the next section, we did not restrict ourselves to only cases
where the hypergeometric series itself converges.  To examine the
effect of the number of coefficients, we considered both $m=30$ and
$m=45$ coefficients.  For each choice of $m$ and $R$, we then chose
1000 sets of seven parameters (and also $z$ if not testing the branch
point) at random as described above.

\begin{table}
  \centering
  \caption{Stability results for recursive computation of the
    asymptotic coefficients at the branch point. All errors are
    geometric means across the sample of the maximum error among the
    $m$ coefficients.}
  \begin{tabular}{|c|cc|cc|} \hline
  & \multicolumn{2}{|c|}{Rel. error in $c_k$} &
\multicolumn{2}{c|}{Rel. error in $\omega_{10}$} \\ \hline
 $R$ & $m =30$ & $m=45$ & $m=30$ & $m=45$ \\ \hline
  1 & $3.9\times 10^{-11}$ & $5.5\times 10^{-11}$  & $2.6\times 10^{-16}$ & $2.6\times 10^{-16}$ \\
  5 & $4.9\times 10^{-12}$ & $1.6\times 10^{-10}$  & $7.5\times 10^{-15}$ & $7.8\times 10^{-15}$ \\
 10 & $3.9\times 10^{-13}$ & $8.7\times 10^{-12}$  & $2.3\times 10^{-13}$ & $1.7\times 10^{-12}$ \\
 50 & $1.9\times 10^{-14}$ & $4.1\times 10^{-14}$  & $1.9\times 10^{-14}$ & $4.1\times 10^{-14}$ \\
100 & $1.8\times 10^{-14}$ & $2.7\times 10^{-14}$  & $1.8\times
10^{-14}$ & $2.7\times 10^{-14}$ \\ \hline
  \end{tabular}
  \label{tab:1}
\end{table}

To quantify the results, for each coefficient we found the relative
error between the value calculated in 53 bit precision and that in 106
bit precision.  We then took the maximum of this relative error out of
all $m$ coefficients.  To obtain some measure of central
tendency among the thousand distinct random samples for a given $m$
and $R$, we took the \emph{geometric} mean of these thousand worst-case
errors, since it is the order of magnitude of the result that is most
important. The results of this calculation are shown in the second and
third columns of table~\ref{tab:1} for the branch point cases (the
unit disk cases are similar and are not shown). We can see already
that there is only a weak dependence of these results on $m$,
indicating that nested recursions are in fact quite stable as we
consider more and more asymptotic coefficients. Somewhat surprisingly,
the average worst-case relative error is seen to be \emph{higher} when
$R$ is smaller, even though as we will see in the next section it is
larger values of $R$ that make the overall computation of the
generalized hypergeometric function more difficult.  We might also be
surprised that the typical worst-case errors can be as large as
$10^{-10}$, since we will see results in the next section that
indicate we can typically calculate the hypergeometric function itself
to higher accuracy in those ranges of $R$.

This is because the relative errors in the asymptotic
coefficients $c_k$ are not by themselves what is most relevant.  More
important is the error estimate $\omega_n$ that we calculate from
those coefficients.  If our error is largest in the highest
coefficients, then that error will be suppressed when we divide by
$n^k$.  Thus, in the fourth and fifth columns we show the relative
error in $\omega_{10}$.  Note that 10 is a very conservative number of
terms to sum; for large parameter values we may typically see many
tens or even hundreds of terms that must be summed to calculate the
hypergeometric function.  Yet we see that we are already very close to
just a couple of orders of magnitude above machine precision for
almost all values of $R$, and in particular the relative errors of
$1\times 10^{-12}$ and $2\times 10^{-14}$ that we shall use in the
next section as relative tolerances to demand of ${}_{q+1}F_q$ are
reasonable. More importantly, the rather remarkable precision with
which we can calculate $\omega_n$ is the key reason that the method of
this paper is much more stable than the $E$-method, as we will see in
section~\ref{sec:e-method}.

\section{Implementation and results}
\label{sec:impl-results}

We now have all of the pieces in place for an acceleration algorithm.
To assemble them, we truncate either equation (\ref{eq:31}) or equation (\ref{eq:28})
for the asymptotics of our partial sums:
\begin{equation}
  \label{eq:43}
  s_n = s + \mu\,\omega_n^{(m)} + O(z^nn^{\lambda-m})
\end{equation}
where we have defined the truncated remainder estimate:
\begin{equation}
  \label{eq:42}
  \omega_n^{(m)} := z^n n^{\lambda} \sum_{k=0}^{m-1} \frac{c_k}{n^k}.
\end{equation}
Here $\lambda$ is $\sigma-1$ away from the branch point, and $\sigma$
at $z=1$. Likewise, the asymptotic coefficients $c_k$ are given by
either (\ref{eq:27}) or (\ref{eq:32}), as appropriate.  We call $m$
the order of our method, and we use equation (\ref{eq:40}) to
calculate the asymptotic coefficients $r_k$ regardless of whether or
not $z=1$. We will need $m+1$ of the $r_k$ to calculate $m$
coefficients $c_k$; though it might seem from
equation~\eqref{eq:32} that we would need $m+2$ of the $r_k$ when at
the branch point, in fact the different appearances of $r_{m+2}$
cancel.

We emphasize that despite our detailed analytic knowledge of
$\omega_n^{(m)}$, it is not by itself a remainder estimate; we must
also know the factor $\mu$.  That cannot be fixed by any of our
analytic calculations; it must be estimated directly from the sequence
of partial sums.  This is a somewhat subtle point, but distinguishes
our method from either Euler-Maclaurin summation \cite{WK03} or
methods based on zeta-functions \cite{boyd09,skoro05,BS06}.  Those
methods each compute an approximation to the remainder directly, and
need the partial sums only to combine with that approximation; our
method requires two computed partial sums so that we may fix $\mu$.
It is in this respect more like fixed-order series acceleration
methods, such as Shanks' method, though we continue to refer to $m$ as
the order for our method, as it is $m$ that determines the rate of
acceleration. 

From equation (\ref{eq:43}) and any two
successive computed partial sums $s_n$ and $s_{n+1}$, we may calculate
an accelerated sum $s^{(m)}_n$ that is hopefully a better
approximation to $s$ than either $s_n$ or $s_{n+1}$:
\begin{equation}
  \label{eq:41}
  s_n^{(m)} = \frac{s_n \omega^{(m)}_{n+1} - s_{n+1} \omega^{(m)}_{n}
  }{\omega^{(m)}_{n+1} - \omega^{(m)}_{n}} .
\end{equation}  
This equation encapsulates our basic method; it is defined only for $m
> 0$. 

Of course, since our estimates for the remainders are only asymptotic,
it may require several computed partial sums before the asymptotic
behavior of the remainder is reached; we expect this number to
increase with increasing $m$, and also to depend on the parameters and
argument of the hypergeometric. Even once we have reached the
asymptotic regime for some particular value of $m$, we may need to
continue calculating estimates for larger $n$ to ensure that the
$O(z^n n^{\lambda-m})$ error term is smaller than whatever tolerance
is desired; again, we expect that the number of such additional
iterations will depend on $z$ and the hypergeometric parameters.
Nevertheless, the basic algorithm contained in equation~(\ref{eq:41})
is ripe for an example, so consider the following evaluation, which we
can do in closed form thanks to Gauss' formula~(\ref{eq:7}):
\begin{equation}
  \label{eq:50}
  {}_{2}F_1\!\left(\left.\begin{matrix}1+4i,1.5+4.5i\\ 
                  3+i\end{matrix}\right|1\right)
 =  -0.003206491294324765 -  0.006293652031968077 i.
\end{equation}
We expect this to be challenging for any series summation technique:
all of the parameters are complex and we evaluate the series at the
branch point; the rate of decay of the remainders is only $n^{-1/2}$.

\begin{figure}
  \includegraphics{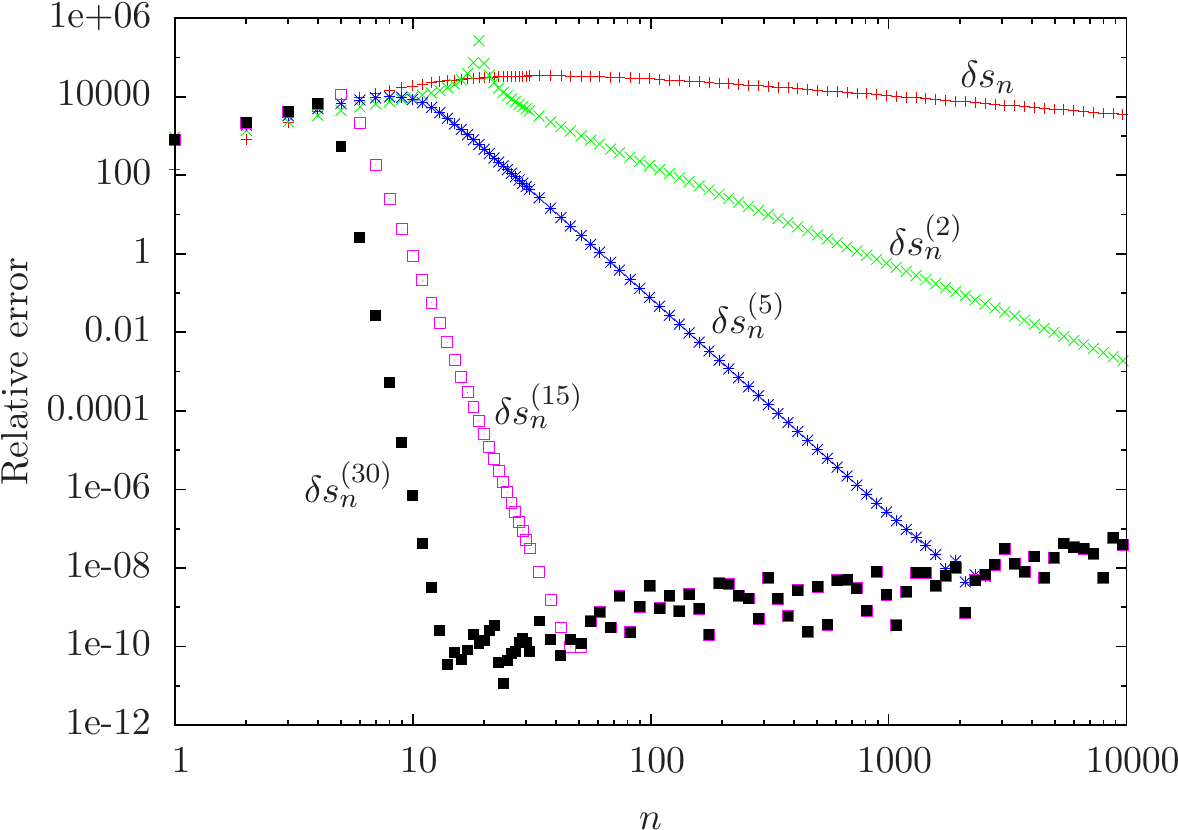}
  \caption{The relative error $\delta s_n^{(m)}$ for different orders
    $m$ of acceleration of the hypergeometric function
    of~(\ref{eq:50}) .}  
  \label{fig:1}
\end{figure}

As we can see in Figure~\ref{fig:1}, however, our method is indeed
able to accelerate the convergence of this hypergeometric function
dramatically. The topmost curve, which shows the relative error
$\delta s_n = |1-s_n/s|$ with no acceleration, has no correct digits
even after summing ten thousand terms, but the acceleration with 30
asymptotic coefficients produces ten digit accuracy with roughly ten
terms.  This figure also shows clearly how the order $m$ of the method
affects the speed of convergence; the linear slopes of the $\delta
s_n^{(m)}$ on a log-log scale illustrate the power law fall-off
predicted by equation~(\ref{eq:31}). Finally, we see that the approach to
convergence need not be monotonic: notice how $\delta s_n^{(2)}$
jumps above the unaccelerated series before it becomes asymptotic.

However, even a small change in the parameters can spoil this
behavior. As a second example, consider the acceleration of
\begin{multline}
  \label{eq:51}
  {}_{2}F_1\!\left(\left.\begin{matrix}1+20i,1.5+25i\\ 
                  3+15i\end{matrix}\right|1\right) \\
    = (-1.508618716765084 + 2.168373234294654 i)\times 10^{-20}.
\end{multline}

\begin{figure}
  \includegraphics{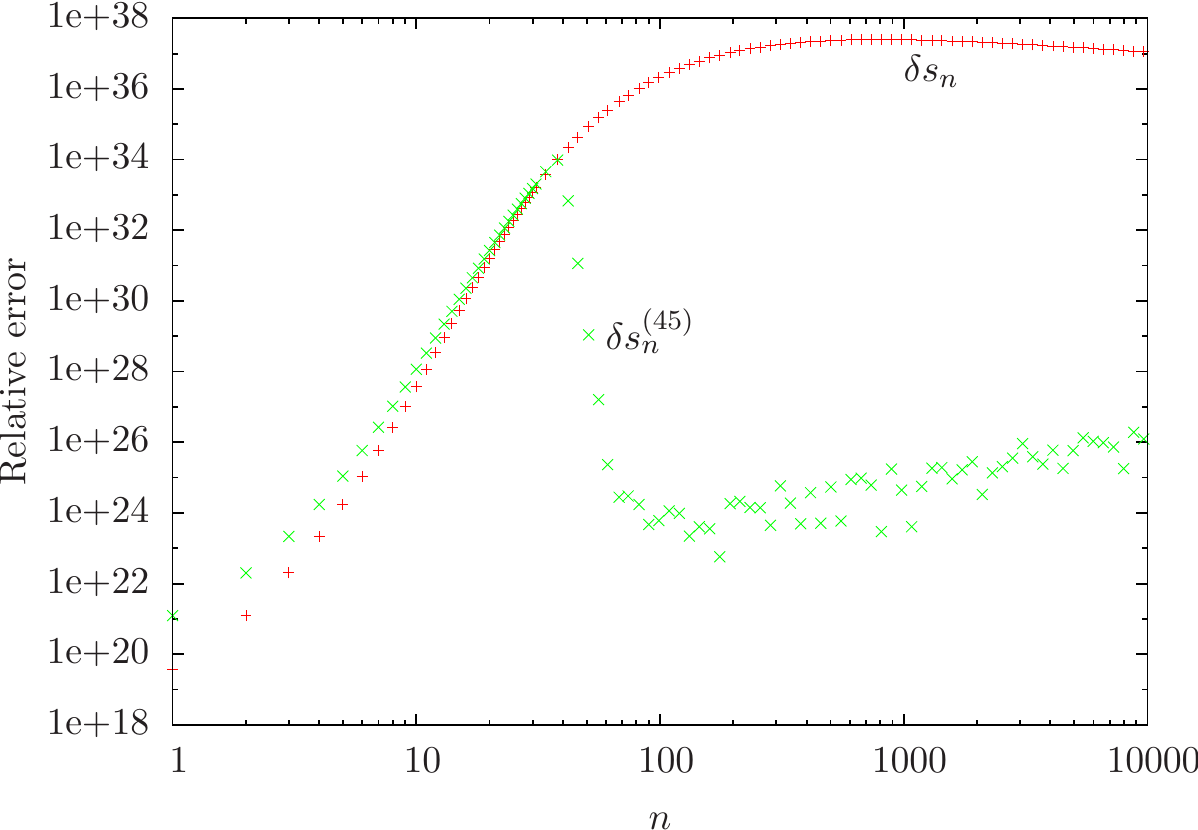}
  \caption{Relative error in the partial sums and accelerated partial
    sums of equation~(\ref{eq:51}). }
  \label{fig:2}
\end{figure}

We plot this in Figure~\ref{fig:2} using 45 asymptotic
coefficients. We do initially see very rapid decrease in the relative
error of the accelerated sum.  However, the relative error then levels
off at about $10^{23}$.  This is despite the very small change made to
the parameters: we have increased the imaginary parts of each
parameter by roughly a single order of magnitude, and we have not
changed the real parts at all.

\subsection{Numerical considerations}
\label{sec:numer-cons}

Of course there is really no mystery here: the largest magnitude that
the partial sum in the evaluation of~(\ref{eq:51}) reaches is roughly
$6.5\times 10^{17}$, or \emph{thirty-seven orders of magnitude greater
than the true value of the function}.  Certainly no method implemented
in fixed IEEE 754 precision can accomplish such an acceleration.

Obviously we could simply implement the algorithm in a symbolic
algebra program that can use arbitrary precision.  However, in at
least some cases this may not be desirable. For instance, the author's
original interest in this project grew out of a problem in
computational quantum gravity that requires the evaluation of
thousands of ${}_4F_3(1)$ with complex parameters as part of a larger
program. Even where higher precision is available it is useful to have
our algorithm correctly determine when that higher precision is really
needed. Finally, because of our precise knowledge of the asymptotics
of the remainder, we can provide an excellent estimate of the error
the algorithm makes when it does converge.

The essential problem is that the accuracy of the algorithm is
determined by its \emph{truncation} error $O(z^n n^{\lambda-m})$ in
equation~(\ref{eq:43}), but before that truncation error becomes
sufficiently small, it may be overwhelmed by the \emph{floating point}
error.  Thus, we need to estimate both sources of error as accurately
as we can, and compare them to decide whether our accelerated
summation has converged to a prescribed accuracy. In the next
two sections we discuss each of these errors in turn. 

\subsubsection{Estimating truncation error}
\label{sec:estim-trunc-error}

In many numerical problems the best estimate of truncation error may
simply be the change between successive estimates, at least when that
change begins to exhibit convergence.  However, because we know
from equation~(\ref{eq:43}) precisely how the truncation error behaves
asymptotically, we can provide a sharper estimate.  Subtracting
successive estimates gives:
\begin{equation}
  \label{eq:47}
  s_{n+1}^{(m)} - s_{n}^{(m)} \approx A z^n n^{\lambda-m}
  \left(z\left(1+\frac{1}{n}\right)^{\lambda-m} -1\right)
\end{equation}
for some constant $A$.  But in fact $A z^n n^{\lambda-m}$ \emph{is} the
leading order of the truncation error, so this suggests that a more 
accurate estimate of the truncation error would be:
\begin{equation}
  \label{eq:48}
  \Delta_{\mathrm{tr}} s_{n}^{(m)} = \frac{s_{n+1}^{(m)} -
    s_{n}^{(m)}}{\left|z\left(1+\frac{1}{n}\right)^{\lambda-m}-1\right|}.
\end{equation}

In our implementation, however, we use the modified estimate:
\begin{equation}
  \label{eq:49}
  \Delta_{\mathrm{tr}} s_{n}^{(m)} = \frac{s_{n+1}^{(m)} -
    s_{n}^{(m)}}{\left|z\left(1+\frac{1}{n}\right)^{-m}-1\right|}.
\end{equation}
Using this estimate greatly reduces the number of false positives
(where the algorithm believes it has converged but has not) when
inside the unit disk $|z|<1$ but with $\Re{\lambda} \gg 0$ (see
section~\ref{sec:conv-as-funct} for details).  Note that
it is only when $|z|<1$ that the real part of the critical exponent
can be positive; at the branch point such series do not converge and
the hypergeometric function is undefined there.  Likewise, even though
the simpler estimate $\Delta_{\mathrm{tr}} s_{n}^{(m)} = 
(s_{n+1}^{(m)}-s_{n}^{(m)})$ has similar rates of convergence or
non-convergence to those we will present in
section~\ref{sec:conv-as-funct}, the estimated error is not nearly as
accurate as with the estimate~(\ref{eq:49}) above, and so we
prefer~(\ref{eq:49}) for that reason also.

\subsubsection{Estimating floating point error}
\label{sec:estim-round-error}

The estimation of floating point error is more complex.  We must be
careful, because overestimation of the floating point error can be
just as dangerous as underestimation: it will cause us to abandon a
calculation that may well have been converging. The very simplest
estimates---for instance, comparing $s_n^{(m)}/\mu$ to the floating
point precision---are disastrous for all but the smallest ranges of
parameters, so we consider a more detailed analysis.  Examining
equation~(\ref{eq:41}), we can identify three primary sources of
potentially significant floating point error in the  calculation of
$s_n^{(m)}$: 
\begin{enumerate}
\item Subtractive cancellation in either the numerator or denominator
  of~(\ref{eq:41}); in each we subtract two nearby floating point
  numbers that we expect will only grow closer to each other as $n$
  increases.
\item Underflow in the calculation of $\omega_n^{(m)}$, particularly if
  $\Re{\lambda} \ll 0$ or $|z| \ll 1$.
\item Accumulated floating point error in $s_n$ and $s_{n+1}$, which
  are calculated recursively using:
  \begin{align}
    s_{n+1} &= s_n + t_n, \label{eq:44a}\\
    t_{n+1} &= z r(n)\, t_n. \label{eq:44b}
  \end{align}
\end{enumerate}

We will need some notation for our analysis. We denote the computed
value of some exact quantity $x$ by $\hat{x}$.  The absolute error
between $x$ and $\hat{x}$ is $\Delta x$, while the relative error is
$\delta x$, so that: 
\begin{equation}
  \label{eq:44}
  \hat{x} = x + \Delta x = x(1+\delta x).
\end{equation}
Applying this to the formula~(\ref{eq:41}) for $s_n^{(m)}$, and
calling the numerator of that formula $\mathcal{N}_n$ and its
denominator $\mathcal{D}_n$, we get:
\begin{equation}
  (\mathcal{D}_n + \Delta \mathcal{D}_n)(s_n^{(m)} + \Delta s_n^{(m)})
  = (\mathcal{N}_n + \Delta \mathcal{N}_n)  \label{eq:45a} \end{equation}
implying:
\begin{equation}
  \Delta s_n^{(m)} = \frac{1}{\hat{\mathcal{D}}_n} (\Delta \mathcal{N}_n - s_n^{(m)}
 \, \Delta \mathcal{D}_n).  \label{eq:45b}
\end{equation}

To derive an estimate from this formula, we must approximate the
errors $\Delta \mathcal{N}_n$ and $\Delta \mathcal{D}_n$.  We estimate
the former by assuming that it is dominated by the error in
calculating the partial sums: 
\begin{equation}
  \label{eq:45}
  \Delta \mathcal{N}_n \approx \Delta s_n |\omega_{n+1}| + \Delta s_{n+1} |\omega_n|.
\end{equation}
We estimate $\Delta \mathcal{D}_n$ in terms of the relative error of the
remainder estimates, rather than absolute.  That is because we wish to
vary our estimate for the remainder error based on whether or not the
calculation underflowed:
\begin{equation}
  \label{eq:18}
  \delta \omega_n = 
  \begin{cases}
    \epsilon_p & \text{if $z^n n^{\lambda}$ is a normalized
      floating point number,} \\
    f_{\mathrm{norm}} |z^{-n} n^{-\lambda}| & \text{otherwise.}
  \end{cases}
\end{equation}
In this equation, $\epsilon_p$ is the machine precision at which
$\omega_n$ is computed, and $f_{\mathrm{norm}}$ is the smallest
normalized number in that precision.  This subtlety is needed because
when $\omega_n$ becomes sufficiently small, it will no longer be
represented by a floating point number with a relative accuracy of
$\epsilon_p$, but rather with a relative accuracy that decreases from
$\epsilon_p$ down to zero.  Such floating point numbers are called
\emph{subnormal}, and $f_{\mathrm{norm}}$ is the smallest floating
point number that is still normalized (possessing the full relative
precision of $\epsilon_p$).

Of course, these estimates involve many approximations.  For example,
the relative error in $\omega_n$ is likely larger than given
by~(\ref{eq:18}) when it is a normalized number, as it is itself
calculated by a non-trivial chain of floating point computations.  But
in those cases it is not a dominant source of error overall, and
so~(\ref{eq:18}) suffices. Combining
equations~(\ref{eq:45b}--\ref{eq:18}), approximating unknown
exact quantities by their computed equivalents when needed, and
taking absolute values to allow for an upper bound on the error, we 
arrive at the final estimate used in our implementation: 
\begin{multline}
  \label{eq:46}
 \Delta_{\mathrm{fp}} s_n^{(m)} = \frac{1}{|\hat{\omega}_{n+1} -\hat{\omega}_n|} \left(
 |\hat{\omega}_{n+1}| \left[\Delta s_n + |\hat{s}_n^{(m)}|\,\delta
 \omega_{n+1}\right] \right. \\
+  \left.|\hat{\omega}_{n}| \left[\Delta s_{n+1} + |\hat{s}_n^{(m)}|\,\delta
 \omega_{n}\right] \right).
\end{multline}

This equation is not complete until we specify how we estimate $\Delta
s_n$. We choose a simple type of  \emph{a posteriori} estimate; for
details, see chapter 3 of  Higham's monograph \cite{higham96}.  The
advantage of such estimates over the more straightforward \emph{a
  priori} estimates is that they take into account cancellation that
occurs during the computation, and are therefore less prone to 
overestimate the error  (though they can be more expensive to compute).

Our estimate starts from equation~(\ref{eq:44a}) and treats
$t_n$ as an exact quantity, accounting only for the accumulation of
error in the recursive sum in $s_n$.  Then it is not hard to
show that the appropriate \emph{a posteriori} estimate starts from
$\Delta s_0 = 0$ and recursively calculates
\begin{equation}
  \label{eq:52}
  \Delta s_{n+1} = \Delta s_n + \epsilon_p\,|s_{n+1}|
\end{equation}
where again $\epsilon_p$ is the machine precision at which the
computation is carried out.

Of course, in reality $t_n$ is also corrupted by ever-growing error,
since it too is calculated recursively.  It is possible to apply a
similar \emph{a posteriori} analysis and derive an estimate for
$\Delta s_n$ that takes this into account. However when this method
was implemented, it tended to severely overestimate the error, greatly
reducing the convergence rate and increasing the number of false
negatives (where the algorithm does not believe it has converged, even
though it has).  At the same time it was more expensive to compute and
even though it could terminate some non-converging cases more rapidly,
its overall performance was slower, even for parameter ranges where a
substantial majority of cases did not converge.  Hence, for the rest
of this paper we consider only the estimate~(\ref{eq:52}) when
estimating the error in the partial sums, and we use that in our
overall estimate~(\ref{eq:46}) for the floating point error in our
accelerated sum.

\subsubsection{Final algorithm}
\label{sec:final-algorithm}

With our estimates for truncation and floating point errors in hand,
we can finally state our complete algorithm, which seeks to
approximate a generalized hypergeometric function to specified
accuracy, or determine that this cannot be done without higher
precision. It takes as input not only the hypergeometric parameters
and argument, but also the order $m$, the desired relative tolerance
$\varepsilon$, and a maximum number of iterations $N$.  The errors are
estimated as described in the preceding two subsections.  In
pseudocode we have:

\begin{tabbing}
\hspace*{2\parindent} \=  \textbf{Input:} $z$, $q$, $\{\alpha_i\}_{i=1}^{q+1}$,
  $\{\beta_i\}_{i=1}^{q}$, $m$, $N$, and $\varepsilon$. \\
\>  \textbf{Calculate:} Exponent $\lambda$ and $\{c_k\}_{k=0}^{m-1}$\\
\>  \textbf{Initialize:} \= $s_0$, $s_1$, and $s_2$; $\omega_1$ and
  $\omega_2$;  $s_1^{(m)}$; $\Delta s_0$\\
\>\textbf{do:}
\> Calculate $s_{n+1}$, $\Delta s_{n+1}$, $r(n)$, and $\omega_{n+1}$ \\
\>\> From $s_{n}$, $s_{n+1}$, $\omega_{n}$ and $\omega_{n+1}$, calculate
$s_{n+1}^{(m)}$ \\
\>\> From $s_{n}$, $s_{n+1}$, $\omega_{n}$, $\omega_{n+1}$, $\Delta s_{n}$,  and
$\Delta s_{n+1}$ calculate $\Delta_{\mathrm{fp}} s_{n}^{(m)}$ \\
\>\> From $s_{n}^{(m)}$ and $s_{n-1}^{(m)}$ calculate
$\Delta_{\mathrm{tr}} s_{n}^{(m)}$ \\
\>\> $n\rightarrow n+1$, $s_{n+1}\rightarrow s_{n}$, $\Delta
s_{n+1}\rightarrow \Delta s_{n}$, and $\omega_{n+1}\rightarrow\omega_n$\\
\>\textbf{while} $n<N$ \textbf{and} $\tau\Delta_{\mathrm{fp}} s_{n}^{(m)} <
\Delta_{\mathrm{tr}} s_{n}^{(m)}$  \textbf{and} $\Delta_{\mathrm{tr}}
s_{n}^{(m)} > \varepsilon\,|s_{n}^{(m)}|$\\
\>\textbf{if} $\Delta_{\mathrm{tr}} s_{n}^{(m)} \le
\varepsilon\,|s_{n}^{(m)}|$ \textbf{return} ``Success'',
$s_{n}^{(m)}$, and $\Delta_{\mathrm{tr}} s_{n}^{(m)}/|s_{n}^{(m)}|$ \\
\>\textbf{else if} $\tau\Delta_{\mathrm{fp}} s_{n}^{(m)} \ge
\Delta_{\mathrm{tr}} s_{n}^{(m)}$ \textbf{return} ``Insufficient
precision'' \\
\>\textbf{else return} ``Maximum iterations reached''
\end{tabbing}

Note in particular that the order of the return statements is such
that if on the same iteration we simultaneously reach a tolerance less
than the specified tolerance but also the floating point error exceeds
the truncation error, we nonetheless deem the acceleration to have
converged.  This is the less conservative approach, but when these
boundary cases are instead treated as failure, we greatly increase the
number of false negatives, without avoiding any additional false
positives. The parameter $\tau$ is an empirical ``fudge factor'' to
lower the estimate of the floating point error; in our final
implementation it was set to 0.1 as this somewhat decreases the false
negative rate. 

\subsection{Testing the method}
\label{sec:testing-method}

The algorithm that we have now derived and described in some detail
was implemented in C, to test its effectiveness as a general purpose
computational strategy for $\fqq$.  Of course we should like to
validate it through testing, but we immediately face the problem of
what to test it against, since the generalized hypergeometric function
can be so challenging to compute in the cases we have in mind.

One choice of course is to test the calculation of $\ftwo$ at the
branch point, since that is simultaneously non-trivial for a series
based computation, yet easily benchmarked against Gauss'
formula~(\ref{eq:7}). That indeed forms the bulk of our test suite, but
we also examined some ${}_3F_2$ and ${}_4F_3$ functions inside the
circle of convergence.  We compared these against the corresponding
calculations from the Python package \textbf{mpmath} \cite{mpmath},
which as mentioned in the introduction takes a fairly sophisticated
approach to calculating generalized hypergeometric functions. However
it was still too slow (at high precision, which we used to ensure
accuracy) to test ${}_3F_2$ and ${}_4F_3$ functions at the branch
point for a large number of randomly generated cases.

All of our test results below will be presented in terms of a
parameter $R$, which sets the scale for choosing the parameters.  The
precise selection was as follows:
\begin{description}
\item[Test cases in the unit disk] We choose a point $z$ with uniform
  probability inside the disk $|z| \le 1$.  We then choose $q+1$ upper
  parameters and $q$ lower parameters, with the real and imaginary
  parts of each chosen uniformly at random between $-R$ and $R$.
\item[Test cases at the branch point] We choose the upper and all but
  one of the lower parameters with both real and imaginary parts
  chosen uniformly in $(-R,R)$, but we choose the last lower parameter
  so that $\Re{\sigma} < 0$.  To do this, if the real part of the sum
  of the upper parameters less the sum of the already chosen lower
  parameters is less than zero, we simply choose the real part of the
  last lower parameter uniformly between that value and $R$.  But if
  that sum is positive, we choose uniformly between 
  that value and either $R$ or that value plus $0.1R$, whichever is
  greater. So the last lower value may have a real part greater than
  $R$.  
\end{description}

With that procedure in mind, there are several questions we wish to
investigate empirically: What is the optimal choice of order $m$? What
is the accuracy of the method, and how often does it terminate
correctly? How accurate is the estimation of error? What are optimal
choices for the maximum number of allowed terms, $N$?  How fast is the
method? We present results on all of these questions in the following
sections; for most of these we only consider ${}_2F_1(1)$.  However
when we examine the overall accuracy and termination we also consider
cases in the unit disk, as described above.

\subsubsection{Convergence as a function of order}
\label{sec:conv-as-order}

We begin with an investigation of the effect of the order on the
convergence of the method.  We have already seen in a simple example
that the order does indeed dramatically affect the rate of
convergence, yet at the same time higher order requires more
precomputation and slows the execution of the method, so we only wish
to invest in this when it is helpful.  But if we choose too low an
order all of the rest of our tests will be essentially meaningless.

Thus, we considered ${}_{q+1}F_q(1)$ for $q = 1$, 2 and 3.  We
generated $10^4$ random cases for each value of $q$ and for each $R$
in $\{1,5,10,50,100\}$.  We required a relative tolerance of
$2\times 10^{-14}$, and then for each randomly chosen parameter set we
called the algorithm with each value of $m$ from 5 to 50.  For each
parameter set for which at least one of these calls converged, we
then observed for which $m$ the fewest number $n$ of partial sums were
needed to achieve the specified convergence; we called this the
optimal order $m_{\mathrm{opt}}$.  We will not present all of the
results here, but simply an illustrative example for the ${}_4F_3$
functions; the results for other choices of $q$ are similar.

\begin{figure}
  \includegraphics{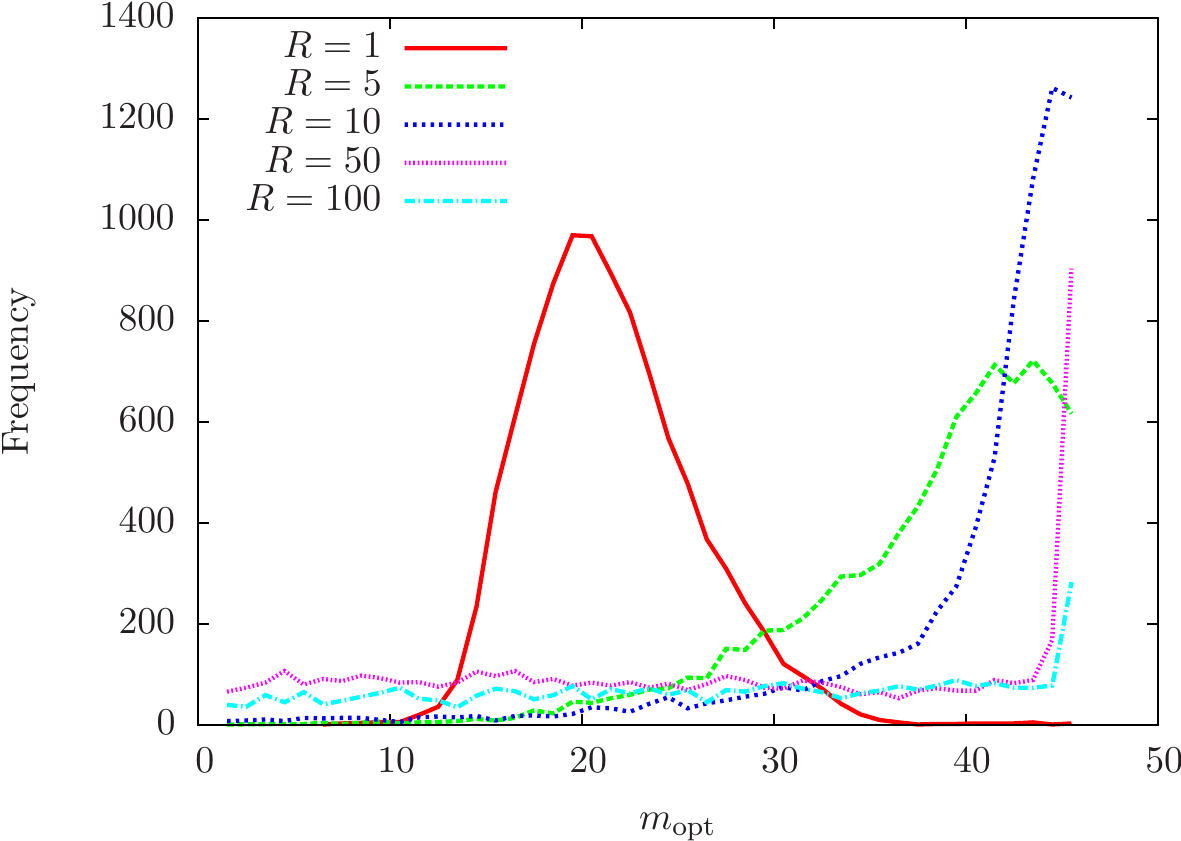}
  \caption{Distribution of optimal order for varying choices of $R$ for the ${}_4F_3(1)$ functions}
  \label{fig:3}
\end{figure}

We can see from figure~\ref{fig:3} that there is a dependence on $R$, but
note that regardless of $R$ the maximum optimal $m$ is 45, even though
we tested up to 50; at these ranges of $R$, at least, there is simply
no point in using more than 45 asymptotic coefficients.  Of course,
this does not prove that using more coefficients will lead to faster
execution, since the increased precomputation may offset the need for
fewer partial sums.  We will examine the speed of the algorithm in
section~\ref{sec:speed-algorithm}. But we do conclude from this plot
and the similar plots for ${}_2F_1$ and ${}_3F_2$ that 45 is a
reasonable upper choice for the order; when we wish to compare the
effect of order on other tests we will also present results
for $m=30$.

\subsubsection{Convergence as a function of the hypergeometric parameters}
\label{sec:conv-as-funct}

By far the most important question for our method is whether or not it
converges, and whether or not it reliably determines when it has. Our
most extensive testing was on this question.

To study it, we again considered $R \in \{1,5,10,50,100\}$, $m=30$ and
$m=45$, and a relative tolerance $\varepsilon$ of either $1\times
10^{-12}$ or $2\times 10^{-14}$.  For each permutation of these
parameter choices, we generated $10^5$ random samples and compared the
computed with the exact value of ${}_2F_1(1)$. We considered four
possible scenarios:
\begin{description}
\item[Convergence] The algorithm claimed that it converged to the
  specified tolerance, and the true relative error was within a factor
  of ten of that specified tolerance.
\item[False positive] The algorithm claimed that it had converged, but
  the true relative error was more than ten times the tolerance.
\item[No Convergence] The algorithm claimed it did not converge, and
  its true relative error was more than the tolerance.
\item[False negative] The algorithm claimed it did not converge, yet
  its true relative error was less than the tolerance.
\end{description}
The factor of ten is somewhat arbitrary but merely reflects the
uncertainty in our truncation error estimate; in this preliminary
investigation it is too stringent to demand that true error be
strictly less than than the tolerance, though we 
will see in section \ref{sec:accur-error-estim} that this is almost
always the case.  Finally, in addition to these four scenarios (which
are mutually exclusive) we also show the percentage of each sample
that reached the maximum number of allowed iterations ($2\times 10^4$
in this set of tests) without converging.  These samples are not
really cases that could have converged had they been given more time,
but rather cases where the floating point error was underestimated and
the algorithm did not terminate early with failure as it should have.

For the case $m=45$, these results are summarized in
tables~\ref{tab:2} and~\ref{tab:3}. We do not present the results for
$m=30$ since they differ very little from these; typically at most a
percent.  We can see from these tables that the convergence rate depends
strongly on $R$, as we would expect from our earlier examples.  But
the false positive rate is extremely low; and even that overstates the
issue: all of the cases labeled as false positives in fact converged
but with a slightly higher ratio between the true and estimated error.
Had we chosen our adjustment factor to be 200 instead of 10, the false
positive rate would be zero in all cases.  Thus, at the branch point
we conclude that when the algorithm terminates with success, it is
essentially always reliable. The false negative rate is also low,
though not nearly so small as the false positive rate.  It is also
somewhat ambiguous, since it tells us only that when the algorithm
terminated with failure, it was in fact really within the prescribed
tolerance; it does not identify scenarios where had the calculation
continued further, convergence would have been achieved. A wide
variation in this rate is therefore possible based on the choice of
floating point error estimate; our implementation uses the choice that
gave the smallest false negative rate of those estimates we
considered. .  Comparing table~\ref{tab:2} to table~\ref{tab:3}, we can
see that requiring higher accuracy does decrease the percentage of
convergent cases, with an effect most pronounced for values of $R$ in
the middle of the range we considered; at high $R$, the percentage of
convergent cases is small enough that most of the effect
of the choice of $\varepsilon$ is masked.

\begin{table}
  \caption{Accuracy of ${}_2F_1(1)$ for $m=45$ and
    $\varepsilon=1\times 10^{-12}$}
  \begin{tabular}{|c|ccccc|}\hline
  $R$ & C & FP & NC & FN & $n_{\mathrm{max}}$ \\ \hline
  1   & 100.0\%& 0.0\%& 0.0\%& 0.0\%& 0.0\%\\ 
  5   & 93.82\%& 0.046\%& 5.91\%& 0.23\%& 5.21\%\\
  10  & 78.89\%& 0.073\%& 20.24\%& 0.80\%& 13.18\%\\
  50  & 36.40\%& 0.014\%& 62.57\%& 1.02\%& 11.90\%\\
  100 & 22.62\%& 0.006\%& 76.56\%& 0.81\%& 9.13\%\\ \hline
  \end{tabular}
  \label{tab:2}
\end{table}

\begin{table}
  \caption{Accuracy of ${}_2F_1(1)$ for $m=45$ and
    $\varepsilon=2\times 10^{-14}$}
  \begin{tabular}{|c|ccccc|}\hline
  $R$ & C & FP & NC & FN & $n_{\mathrm{max}}$ \\  \hline
  1   & 98.90\%& 0.012\%& 1.08\%&  0.01\%& 0.90\% \\
  5   & 81.39\%& 0.096\%& 17.97\%& 0.54\%& 14.64\% \\
  10  & 65.90\%& 0.095\%& 32.87\%& 1.13\%& 19.59\% \\
  50  & 29.96\%& 0.019\%& 68.66\%& 1.36\%& 12.26\% \\
  100 & 18.53\%& 0.009\%& 80.70\%& 0.77\%& 9.20\% \\ \hline
  \end{tabular}
  \label{tab:3}
\end{table}

We also investigated the accuracy for points chosen in the unit disk.
Here for simplicity we considered only $m=45$ and $\varepsilon=2\times
10^{-14}$, but we were able to examine ${}_3F_2$ and ${}_4F_3$
functions as well, though for smaller ranges in $R$.  These results
are shown in table~\ref{tab:4}.  We see that as we might expect the
overall convergence rates are higher than for the corresponding value
of $R$ at the branch point, and though the data is somewhat limited
there does not seem to be a strong dependence on the order $q$ of the
hypergeometric; certainly not nearly as strong as the dependence on
$R$. One subtlety not shown in this table is that it is no longer true
that the (albeit rare) false positives are necessarily benign; roughly
half of the false positives for the $R=100$ cases of ${}_2F_1$, for
instance, had fewer than half of the digits correct; in many cases no
digits correct.  This always happens when the exponent
$\lambda$ has a large positive real part; as already mentioned in
section~\ref{sec:estim-trunc-error}, our choice of estimate for the
truncation error reduces the fraction of these false positives.  We
can now quantify that assertion: had we used the truncation error
estimate~(\ref{eq:48}) instead of~(\ref{eq:49}), then even with an
adjustment factor of one thousand instead of ten, our percentage of
false positives for the $|z|<1$, $R=100$ case of ${}_2F_1$ would be 
11.22\%, or one hundred times greater.  The false positives that 
still persist even using~(\ref{eq:49}) are those where the terms of
the series temporarily become much smaller (by many orders of
magnitude) before again increasing.  This causes the series to appear
to be rapidly converging when it is not, and our acceleration method
is unable to distinguish this from true convergence.  But neither are
other methods; even commercial algebraic packages were found to
falsely return convergence on these cases, unless specifically
instructed to calculate results to very high precision.  We know of no
reliable way of deciding that this will happen; we must just
fortuitously choose to calculate at sufficiently high precision.   

\begin{table}
  \caption{Accuracy of ${}_{q+1}F_q(z)$ calculations for $|z|<1$ with
    $m=45$ and $\varepsilon=2\times 10^{-14}$.} 
  \begin{tabular}{|c|c|ccccc|}\hline
   & $R$ &  C & FP & NC & FN & $n_{\mathrm{max}}$ \\  \hline
  \multirow{5}{*}{${}_2F_1$} &   1 & 99.86\%& 0.002\%& 0.14\%& 0.0\%& 0.0\% \\
                             &   5 & 94.63\%& 0.15\%& 5.10\%& 0.12\%& 0.011\% \\
                             &  10 & 85.98\%& 0.13\%& 13.70\%& 0.19\%& 0.036\% \\
                             &  50 & 49.32\%& 0.12\%& 50.36\%& 0.20\%& 0.17\% \\
                             & 100 & 31.73\%& 0.11\%& 68.03\%& 0.13\%&
                             0.11\% \\ \hline
  \multirow{3}{*}{${}_3F_2$} &  1 &  99.76\%& 0.03\%& 0.21\%& 0.0\%& 0.0\% \\
                             &  5 & 92.85\%& 0.20\%& 6.78\%& 0.10\%& 0.01\% \\
                             & 10 & 84.88\%& 0.13\%& 14.83\%& 0.16\%& 0.03\% \\ \hline
  \multirow{2}{*}{${}_4F_3$} &  1 & 99.74\%& 0.0\%& 0.25\%& 0.010\%& 0.0\% \\
                             &  5 & 91.35\%& 0.12\%& 8.44\%& 0.09\%& 0.03\% \\ \hline
  \end{tabular}
  \label{tab:4}
\end{table}

In summary, except for those very few cases just mentioned, the method
is highly effective at either accelerating the series or determining
that a higher precision is needed; the error estimate it returns is
almost always accurate to within a factor of ten.  If the arguments to
the function are themselves much larger than about ten, then it is
increasingly unlikely that the method will converge in double
precision, but it will correctly identify that failure, and a
sufficiently high precision implementation should succeed.

\subsubsection{Accuracy of the error estimate}
\label{sec:accur-error-estim}

\begin{figure}
  \centering
  \includegraphics{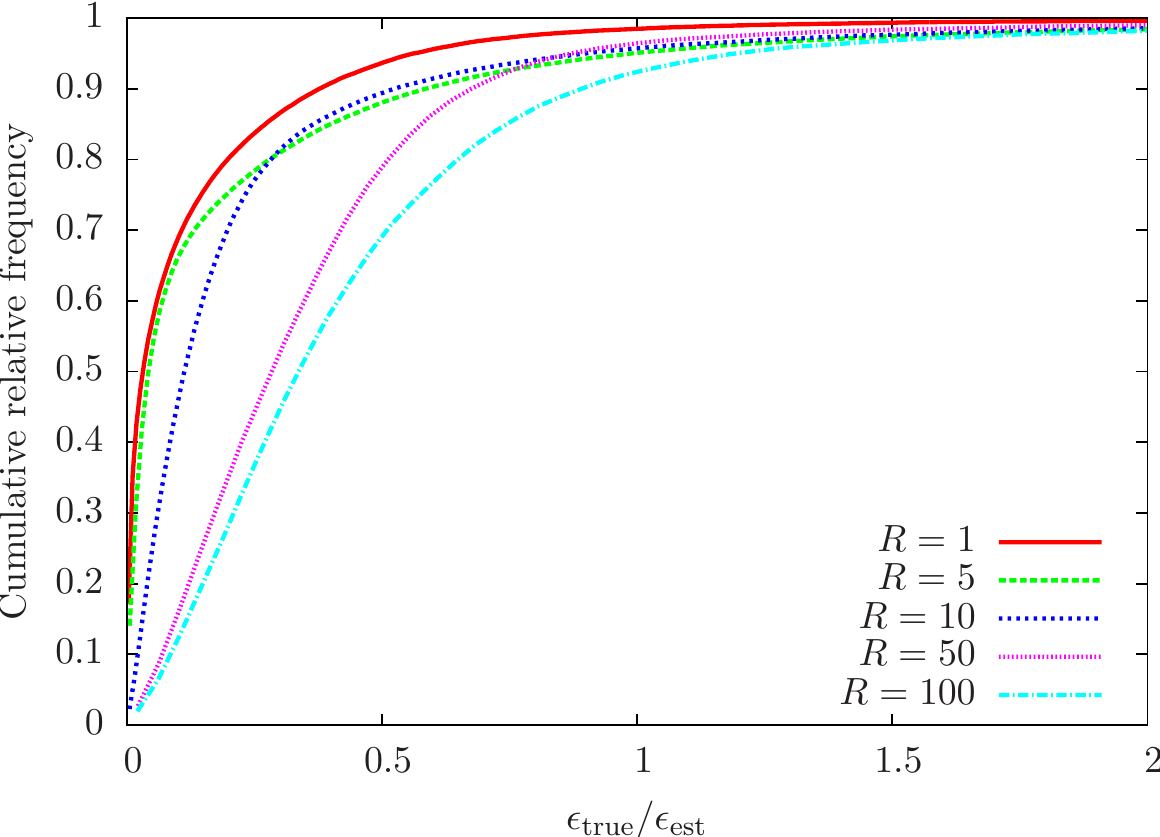}
  \caption{The cumulative relative frequency of the ratio of the true
    error $\epsilon_{\mathrm{true}}$ to the estimated error among
    those of the randomly chosen calculations that converged, when the
    order $m$ is 30 and the requested relative tolerance is $2\times
    10^{-14}$.}  
  \label{fig:4}
\end{figure}

The results of the previous section clearly indicate that our
truncation error estimate is quite reliable, but it is useful to
consider this in more depth.  Thus, from each of the trials that
converged or were false positives in the tests of the previous
section, we can examine the ratio of the true relative error to that
estimated by the algorithm from equation~(\ref{eq:49}).  As an example
of this behavior, we show in figure~\ref{fig:4} a relative cumulative
frequency plot for each of the five values of $R$ we tested, for the
${}_2F_1(1)$ function with $m=30$ and $\varepsilon=2\times 10^{-14}$.
As this figure shows, the error estimate is excellent and only weakly
depends on $R$; in 90\% or more of cases the estimated error is less
than the true error, and for essentially all cases it is within a
factor of two.

\subsubsection{Number of partial sums needed}
\label{sec:number-partial-sums}

\begin{figure}
    \includegraphics{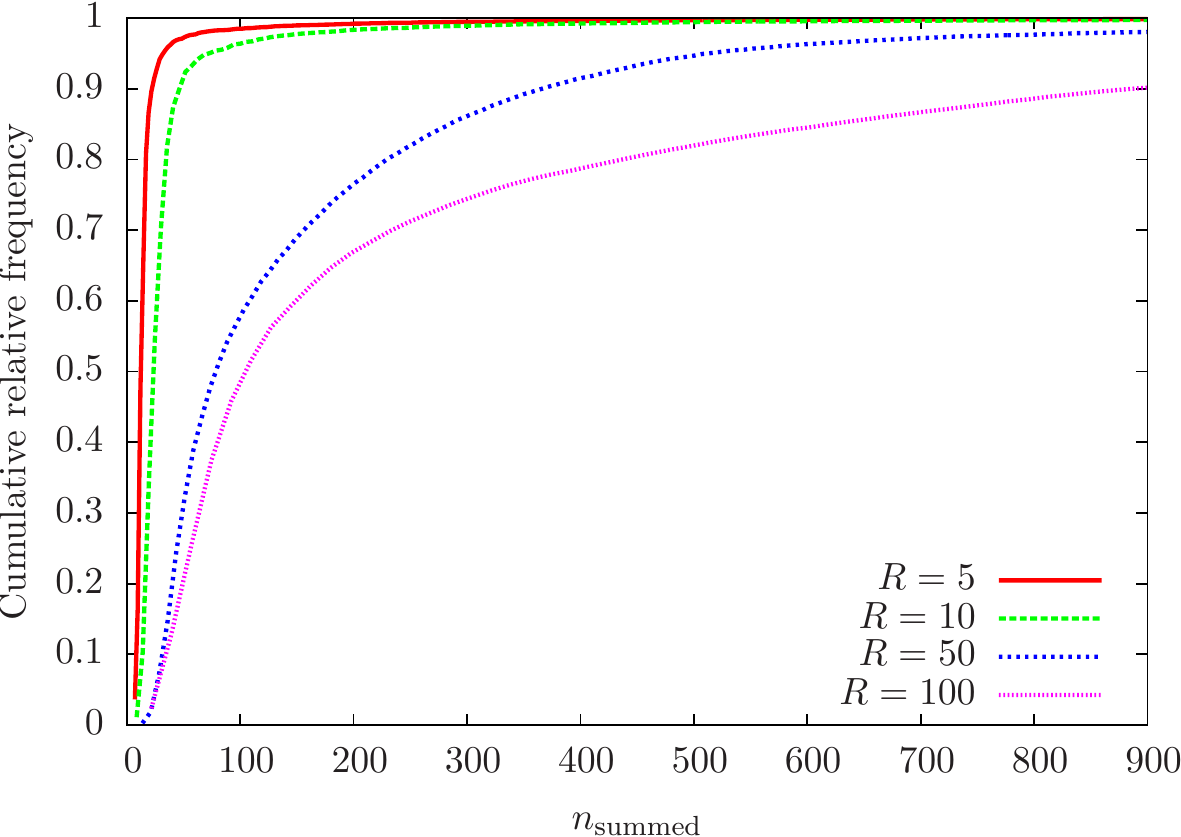}
  \caption{The cumulative relative frequency of the number of terms
    $n$ that were needed among those of the randomly chosen
    calculations that converged, when the order $m$ is 30 and the
    requested relative tolerance is $2\times 10^{-14}$.} 
  \label{fig:5}
\end{figure}

Another important parameter of the method is the maximum number of
iterations allowed before the method is deemed to have failed,
regardless of the error estimate.  Since up to 20\% of cases may reach
this limit, it is important not to have it unnecessarily high, as
otherwise we waste time on an unsuccessful calculation.  A sample of
such behavior---again as a cumulative relative frequency plot---is
shown in figure~\ref{fig:5}; results for other choices of $m$ or
$\varepsilon$ are similar.  Unlike the error ratio, we see a much
stronger dependence on the size of the parameters, as measured by $R$,
but we can see that at least for $R$ up to 100 choosing the maximum
$N$ to be about one thousand is quite conservative; for smaller
parameter ranges this can be reduced even further.

\subsubsection{Speed of the algorithm}
\label{sec:speed-algorithm}

Finally, apart from the accuracy of the algorithm its most essential
characteristic will be its speed. To examine this, we generated one
thousand sets of parameters for each of our usual set of five $R$ values,
and then averaged the time of each parameter set from one hundred runs
on a 1000 MHz AMD Athlon processor under Linux.  We report separately
the time for all runs versus just those runs that converged, and we
set the relative tolerance at $2\times 10^{-14}$ and the maximum
number of allowed iterations at $2\times 10^3$.  For evaluations at
the branch point, we obtain the results shown in table~\ref{tab:5}.

\begin{table}
  \caption{Speed of the algorithm at the branch point, when $N=2\times
    10^3$ and $\varepsilon=2\times 10^{-14}$.}
  \begin{tabular}{|c|c|cc|cc|} \hline
\multicolumn{2}{|c|}{} & \multicolumn{2}{c|}{$m=30$} &
\multicolumn{2}{c|}{$m=45$} \\ \hline
    Function & $R$ & All & Only converged & All & Only converged \\ \hline
   \multirow{5}{*}{${}_2F_1(1)$} &   1 & 0.154 ms & 0.0951 ms & 0.307 ms & 0.250 ms \\
                                 &   5 & 0.675 ms & 0.115 ms  & 0.858 ms & 0.292 ms \\
                                 &  10 & 0.911 ms & 0.145 ms  & 1.15 ms  & 0.339 ms \\
                                 &  50 & 0.982 ms & 0.332 ms  & 1.25 ms  & 0.505 ms \\
                                 & 100 & 1.11 ms  & 0.539 ms  & 1.38 ms  & 0.728 ms \\ \hline
   \multirow{5}{*}{${}_3F_2(1)$} &   1 & 0.119 ms & 0.109 ms  & 0.338 ms & 0.330 ms \\
                                 &   5 & 0.680 ms & 0.124 ms  & 0.843 ms & 0.307 ms \\
                                 &  10 & 0.966 ms & 0.157 ms  & 1.19 ms  & 0.323 ms \\
                                 &  50 & 1.15 ms  & 0.330 ms  & 1.39 ms  & 0.485 ms \\
                                 & 100 & 1.21 ms  & 0.406 ms  & 1.46 ms  & 0.546 ms \\ \hline
   \multirow{5}{*}{${}_4F_3(1)$} &   1 & 0.118 ms & 0.0981 ms & 0.322 ms & 0.310 ms \\
                                 &   5 & 0.726 ms & 0.118 ms  & 0.855 ms & 0.258 ms \\
                                 &  10 & 1.05 ms  & 0.148 ms  & 1.29 ms  & 0.308 ms \\
                                 &  50 & 1.21 ms  & 0.286 ms  & 1.46 ms  & 0.438 ms \\
                                 & 100 & 1.33 ms  & 0.381 ms  & 1.58 ms  & 0.546 ms \\ \hline
  \end{tabular}
  \label{tab:5}
\end{table}

We see from table~\ref{tab:5} that as we would
expect, increasing the order of the algorithm does increase the
execution speed, but more so for the smaller parameter
choices, where the algorithm converges quite quickly and there is
little to be gained by using more asymptotic coefficients.  Within a
given order $m$, we also see the time increase with $R$, by roughly a
factor of two as we move from $R=5$ to $R=100$.  The results for cases
inside the unit disk are similar, except that the increase in running
time as we increase $m$ is more significant; at least a factor of
two for all values of $R$, not just small values.

\section{Comparison with other methods}
\label{sec:comp-with-other}

While the tests of the previous section are a convincing validation of
the method of this paper, it is also important to compare the method
to others.  Here we are somewhat hampered: most of the literature on
computing the generalized hypergeometric function gives only a few
examples, and not the kind of large, randomly selected test cases we
used for testing. Also interest is often focused on obtaining the
greatest accuracy with the fewest terms summed, and robustness
considerations (such as automatic termination) are rarely mentioned.
For these reasons, in this section we use a different implementation
of our algorithm, in Python, using mpmath \cite{mpmath} to provide
arbitrary precision arithmetic.  By using an interpreted language and
software implementations of higher than machine precision, we of
course pay a large performance penalty (typically a factor of roughly
a thousand).  However most of the other
methods we consider here are also implemented in such systems, so a
comparison between the two is still meaningful.

Moreover, our interest here is centered on methods that can handle
generic, complex parameters.  There can certainly be
particular choices of parameters (and argument) for which other methods
are more efficient, but when all of the parameters and the argument
are permitted to be complex, the number of possible ``special cases''
grows bewilderingly large.  So we focus on two methods that are
proposed as general-purpose algorithms for hypergeometric functions,
and also on the $E$-method (mentioned in section
\ref{sec:surv-series-accel}) which bears superficial similarity to the
method of this paper.

\subsection{Zeta function acceleration}
\label{sec:zeta-funct-accel}

This method \cite{skoro05,BS06} is designed to evaluate
${}_{q+1}F_q(1)$, and as such is perhaps still something of a
special-purpose algorithm.  But the branch point is the most
challenging case, and the authors consider complex parameters as well
as several optimizations of their method.  That method is based on
directly summing the first $N$ terms of the series, and then
approximating the remainder in terms of $m$  Hurwitz zeta functions.
The method of \cite{BS06} extends that of \cite{skoro05} by allowing
for a complex parameter $\alpha$ that is determined through a symbolic
algebra problem requiring the solution of a nonlinear optimization
through Gr\"{o}bner bases. The complexity of this optimization problem
grows with $m$, so the authors do not consider values of $m$ as large
as those we can easily handle with our method. 

In the optimized work \cite{BS06}, a few examples are considered,
and here we compare two of them with our method. The
timing results of \cite{BS06} were for an AMD Athlon64 3500+
processor; our processor is comparable if a little faster (it is a
4600+ model).

The first example considered in \cite{BS06} are ${}_3F_2(1.6+7i,
2.4-i,\sqrt{2};3+i,\sqrt{6}+i;1)$. The method of that paper is able to
evaluate that function to ten digit accuracy with $N=35$ in 0.3
seconds; our method required $N=17$ and 0.37 seconds.  But as the
required precision is increased, the advantage of our method grows: 15
digit accuracy with the zeta-function method required $N=100$ and
1.1 seconds, but for our method $N=25$ and still 0.37 seconds; 35
digit accuracy required $N=3500$ and 14 seconds for them, but $N=110$
and still just 0.37 seconds for us.

In fact in every case considered in \cite{BS06} our method
out-performed that method, often substantially; to keep our discussion
brief we consider just one more example.  The most challenging example
considered in that paper was ${}_4F_3(2.4+30i,-0.3+0.5i,2.2-i,0.5+i;
1.8, 1.1-i,2+17i;1)$.  The authors could achieve 10 digit accuracy
with roughly 1000 terms summed, whereas we achieve the same
with only 136 terms.  To achieve 20 digit accuracy they required
approximately 6000 terms, whereas we needed only 290.  It is true that
we can use a larger value of $m$ than those authors (we used 30 in our
tests; they used either 10 or 15) but that is again because it is easy
for us to increase $m$, as there is no system of $m$ polynomials to
solve in our method.

\subsection{Euler-Maclaurin summation}
\label{sec:euler-macl-summ}

Euler-Maclaurin summation is based on a specific analytic form of the
remainder of a series, expressed as an integral and a weighted sum of
derivatives of the terms \emph{with respect to the term index}
\cite{WK03}: 
\begin{equation}
  \label{eq:53}
  \rho_n = \int_{n+1}^{\infty} t(k)\,dk + \frac{1}{2}t(n+1) -
  \sum_{j=1}^{m} \frac{B_{2j}}{(2j)!}t^{(2j-1)}(n+1).
\end{equation}
Here the $B_{2j}$ are the Bernoulli numbers, and we have assumed that
not only do the terms $t_k$ go to zero as $k\rightarrow\infty$, but so
also do all of the derivatives of the terms with respect to the term
index; the method is easily generalized when that does not hold.

It is this need to integrate and differentiate terms with respect to
the term index that makes this method challenging.  For the Riemann
zeta function, the authors of \cite{WK03} could carry this out
analytically, but for generalized hypergeometric functions an analytic
solution seems intractable.  However a numerical implementation of
this method underlies the \textbf{mpmath} \cite{mpmath} calculation of
generalized hypergeometric functions near the branch point, so we make
our comparison with that implementation.

The first test case considered in \cite{mpmath} is
${}_4F_3(\frac{1}{3},1,\frac{3}{2},2;\frac{1}{5},
\frac{11}{6},\frac{41}{8};1)$. For mpmath and 25 digit accuracy, its
Euler-Maclaurin based summation method requires 2.4 seconds, while on
the same machine our method requires only 0.5 seconds.  But if we try
to extend the test cases of the previous subsection, then the
Euler-Maclaurin based approach is completely incapable of grappling:
the first and simplest test case we considered runs for several
minutes before returning a failure to converge.  It fares even worse
in the other test cases.

\subsection{$E$-method}
\label{sec:e-method}

Unlike the method of this paper, the previous two methods require the
computed partial sum only so they can add their estimates of the
remainder to that partial sum; the remainder itself they calculate
without direct reference to the sequence of partial sums.  Our method
requires two successive partial sums, because both $s$ and $\mu$ in
equation \eqref{eq:43} are unknown; we are effectively solving a
$2\times 2$ linear equation.  At the other extreme, we could use only
our knowledge of the leading behavior of the remainder, and rather
than precomputing the asymptotic coefficients $c_k$, we can determine
the coefficients $\tilde{c}_k$ in:
\begin{equation}
  \label{eq:54}
  s_n = s + \sum_{k=1}^{m} \tilde{c}_k \frac{z^n n^{\lambda}}{n^{i-1}}.
\end{equation}
Here the coefficients $\tilde{c}_k$ are related to our asymptotic
coefficients $c_k$ through $\tilde{c}_k = \mu c_{k-1}$.  

This approach is the $E$-method already mentioned in the introduction,
for the particular choice of functions $g_k(n) = z^n /
n^{k-\lambda-1}$. That method is described in the monograph
\cite{BR91} of Brezinski and Redivo Zaglia and was independently
discovered by Schneider \cite{schneider75}, H\r{a}vie \cite{havie79}
and Brezinski \cite{brezinski80}; a stable numerical implementation is
described by Brezinski in \cite{brezinski82}.  

In fact, the $E$-method is perhaps more properly thought of as a class
of methods; most existing methods can be subsumed by specializing to a
particular choice of the $g_k(n)$. Indeed Levin's original work
\cite{levin73} on nonlinear sequence transformations can be analyzed
as a specialization of the $E$-method to an asymptotic expansion in
inverse powers, with differing simple remainder estimates that enable
application to a variety of different sequences.  The algorithm
we analyze now is another such specialization, where the remainder
estimate is given by our analytic knowledge of only the leading order
of the asymptotic truncation error.  Thus the comparisons of this
section can also be considered a comparison to a variant of Levin's
methods.  

Not only is such a specification necessary to have a complete
algorithm, but particular specializations will often allow simpler
implementations than those described in the references above for the
general-purpose $E$-method.  In our present case, the special form of
our $g_k$ makes it simpler to use the recursive scheme in section 7.2
of Weniger's review article \cite{weniger89}, itself based on the work
of Fessler \textit{et al.} in \cite{FFS83}.   

We compare first the complexity of the two approaches.  If we
assume that $m \gg q$, then the complexity of the algorithm of this
paper is roughly $\frac{1}{3}m^3+10m^2+(2m+7)N$, while for this
implementation of the $E$-method it is $6m^2 N$. Thus, as $N$ grows
beyond $m$, the method of this paper has clearly better complexity;
even when $m\approx N$ it is somewhat superior.  Note that the
$E$-method cannot have $N<m$, since we must always consider at minimum
$m+1$ partial sums.

\begin{figure}
  \centering
  \includegraphics{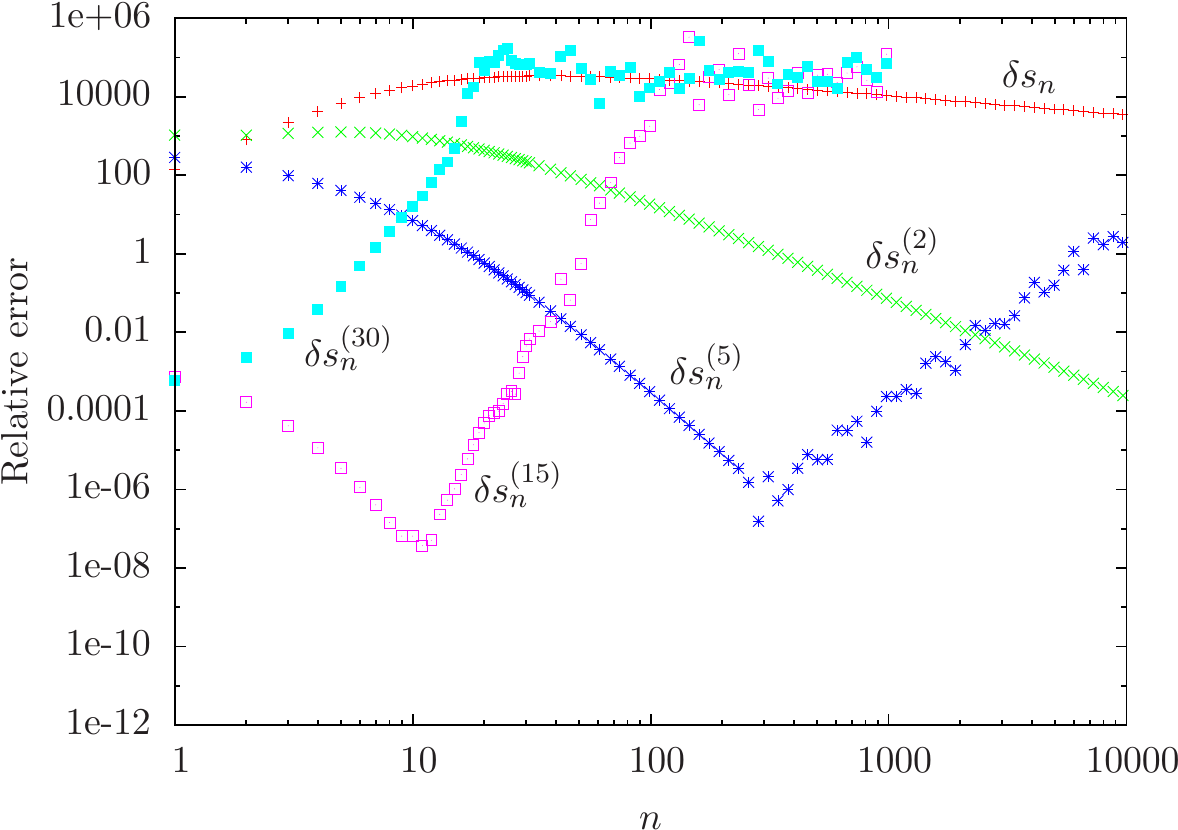}  
  \caption{The relative error $\delta s_n^{(m)}$ for different orders
    $m$ of the $E$-method applied to (\ref{eq:50}) .}  
  \label{fig:6}
\end{figure}

The real advantage of our method, however, is in its stability.  To
illustrate that, we return to the first example we considered, from
equation \eqref{eq:50} and as shown in figure \ref{fig:1}.  Figure
\ref{fig:6} shows the corresponding plot for the $E$-method using
calculations in 80 bit (long double) precision.
Comparing to figure \ref{fig:1}, we note several differences:
\begin{enumerate}
\item The overall error is larger.
\item The floating point error grows much more quickly with the order
  $m$ of the transformation.
\item Once the minimum error is reached, the error rapidly begins
  increasing again, so that automatic termination would be much more
  challenging to implement.
\end{enumerate}

We can understand this overall loss of stability if we consider that
the $E$-method is essentially solving the linear system:
\begin{multline}
  \label{eq:55}
  \begin{bmatrix}
    1 & z^n/n^{-\lambda} & z^n/n^{1-\lambda} & \ldots &
    z^{n}/n^{m-1-\lambda} \\
    1 &  z^{n+1}/(n+1)^{-\lambda} & z^{n+1}/(n+1)^{1-\lambda} & \ldots
    & z^{n+1}/(n+1)^{m-1-\lambda} \\
   \vdots & \vdots & \vdots & \ddots & \vdots \\
    1 &  z^{n+m}/(n+m)^{-\lambda} & z^{n+m}/(n+m)^{1-\lambda} & \ldots
    & z^{n+m}/(n+m)^{m-1-\lambda} 
  \end{bmatrix} \\
  \times 
  \begin{bmatrix}
    s \\ \tilde{c}_1 \\ \vdots \\ \tilde{c}_m
  \end{bmatrix} 
 =
 \begin{bmatrix}
   s_n \\ s_{n+1} \\ \vdots \\ s_{n+m}
 \end{bmatrix}.
\end{multline}
Of course, this system is not explicitly solved at each step, since it
is only $s$ that we need, but it is still the stability of the
underlying system \eqref{eq:55} that dictates the stability of the
recursive scheme for $s$. We can quantify that through condition
numbers. Examining figure \ref{fig:1}, we see that we would expect to
achieve 10 digit accuracy when $n=10$ if our order is $m=30$.  Yet the
condition number of the matrix of \eqref{eq:55} for those choices is
$1.3\times 10^{28}$, far too large to allow a solution in long double
precision.  If we try to avoid this by decreasing $m$, then we must
also increase $n$; again from figure \ref{fig:1} we estimate that if
$m=15$ we would need roughly $n=40$; now the condition number is
$3\times 10^{89}$.

These large condition numbers are not coincidental.  We saw that our
method begins to lose precision whenever $|s| \ll |\mu|$, since then
the unknown that we care about in our linear system becomes much
smaller than the other unknown; this must spring from instability in
the underlying ($2\times 2$) system we are solving.  With the
$E$-method, the same problem can arise if $|s| \ll \tilde{c}_k$ for
any of the $c_k$ that we solve for.  But that will happen generically:
the coefficients are asymptotic and grow rapidly with $m$; for the example
we consider here we have $|c_{30}| = 10^{18}$.  Hence the
corresponding linear systems must be unstable.

From this perspective, the chief advantage of the method of this paper
is that it bypasses such an unstable linear system.  Instead, as we
saw in section \ref{sec:stability-recursions}, we can determine the
$\omega_n$ to high accuracy, and that accuracy does not decrease
rapidly with $m$, and actually increases with $n$.  For this reason,
our method is much more stable.

\section{Conclusions}
\label{sec:conclusions}

Summarizing, we have shown that it is possible to derive the
asymptotics of the remainders of the partial sums of the generalized
hypergeometric function ${}_{q+1}F_q$ to any desired order in inverse powers of
$n$. We have given explicit formulas for the remainders in terms of
the hypergeometric parameters and argument. This analytic result is the
basis for a new series acceleration technique that can dramatically
accelerate the convergence of the generalized hypergeometric series,
making it feasible to evaluate these for complex arguments, even at
the branch point $z=1$.  As implemented in C, the algorithm can be
limited by the fixed precision of standard floating point types, but
even in this case the precise asymptotic knowledge available enables
us to determine correctly when the acceleration has converged.  The
method seems much more efficient and robust than any others we have
found in the literature that are applicable to ${}_{q+1}F_q$ at
generic complex arguments and parameter.

There are still some open issues, which are natural starting points
for future research:
\begin{itemize}
\item At present, the algorithm is very slow \emph{near} the branch
  point, much more so than \emph{at} the branch point.  As shown by
  B\"{u}hring \cite{buhring03,buhring92,buhring87,BS98} and N{\o}rlund
  \cite{norlund55}, there is a close connection
  between the asymptotics of the partial sums at the branch point and
  the behavior of the function near the branch point. It would be
  interesting to see if this can be leveraged to evaluate the function
  near the branch point using the (faster) evaluation at the branch
  point; of course this is not just a simple series expansion
  precisely because we are near a singularity.
\item As we have noted throughout this paper, our C implementation is
  limited by fixed floating point precision.  Of course it is
  straightforward to implement the algorithm in any of the free or
  commercial symbolic computational programs that support arbitrary
  precision, but it could also be useful to continue development of an
  arbitrary precision routine in a low-level language, by taking
  advantage of existing higher precision libraries like MPFR
  \cite{FHLPZ07} or QD \cite{HLB01}.  Moreover, the acceleration
  presented here is almost certainly not ideal for all inputs; for
  small $|z|$, for instance, there is no need to use any acceleration
  at all.  Even when the method of this paper is best, a more
  automatic implementation should choose the order $m$ and maximum
  iterations $N$ to minimize the computation required to achieve the
  desired accuracy.  All of these considerations together could lead
  to a reliable and fast library for generalized hypergeometric
  evaluations. 
\item We have focused in this paper on the case $p=q+1$ because that
  restriction was needed to apply the results of \cite{WL92a}.
  However the general case can be handled by including the further
  results of those same authors in \cite{WL92b}; work on this
  extension is already underway. 
\item We have also only described an algorithm in which we use the
  precise remainder estimates and a constant correction factor $\mu$,
  in contrast to traditional series acceleration techniques that use
  simple remainder estimates and sophisticated correction factors.
  But the two choices are not exclusive, and it would be interesting
  to investigate the performance of a method that combines the
  remainder estimates of this paper with the correction factors
  $\mu_n$ of traditional series acceleration techniques.
\item We have limited ourselves to asymptotic expansions in
  inverse powers of $n$, because those are the asymptotic functions in
  which the results of \cite{WL92a} are couched.  However, Weniger has
  found \cite{weniger07,weniger10} that inverse factorial series can
  also be very useful---in some cases much more powerful---than
  inverse powers, and it is worth investigating if a similar expansion
  would be effective here. In particular, it is shown in
  \cite{weniger07,weniger10} than an expansion in inverse powers can
  be transformed to an expansion in inverse factorial series, so the
  question is really how the stability and efficiency of such a scheme
  compares with the method presented here.
\item Finally, it is worth investigating how well the method of this
  paper performs when applied to other functions for which an analytic
  asymptotic expansion of the term ratio is available.  Though many
  series of practical interest in the sciences are available only
  numerically or as expensive computations (e.g., perturbation
  series), there are still other series of practical interest where we
  have available the necessary analytic knowledge, and would like to
  take advantage of that knowledge to efficiently and robustly
  evaluate the functions.  
\end{itemize}

\section*{Acknowledgements}
This research was supported in part by the Math/Science research fund
of Abilene Christian University, and that support is gratefully
acknowledged. I also thank the Max-Planck-Institut f\"{u}r
Gravitationsphysik, Hannover, Germany, for their hospitality and
support while this paper was written.  I also thank Rafa{\l} Nowak for
pointing out a sign error in an earlier draft of this paper, and the
two anonymous referees for their several suggestions, including in
particular the suggestion of several references, and other methods
with which to compare this work.

\bibliographystyle{unsrtnat}
\bibliography{jlw_ghg_accel}

\end{document}